\documentclass[11pt,reqno]{amsart}
\usepackage{amssymb,mathrsfs,graphicx,subfigure, enumerate}
\usepackage{amsmath,amsfonts,amssymb,amscd,amsthm,bbm}
\usepackage{extpfeil}
\usepackage{graphicx,colortbl}
\usepackage{float}
\usepackage{epsfig}
\usepackage{caption}
\usepackage{bm}
\graphicspath{{Figure/}}

\usepackage[margin=1in]{geometry}

\makeatletter
\@namedef{subjclassname@2020}{\textup{2020} Mathematics Subject Classification}
\makeatother

\title[Navier-Stokes equations in half space]{Asymptotic behavior toward viscous shock for impermeable wall and inflow problems of barotropic Navier-Stokes equations}

\author[Huang]{Xushan Huang}
\address[Xushan Huang]{\newline Department of Mathematical Sciences \newline Korea Advanced Institute of Science and Technology, Daejeon  34141, Republic of Korea}
\email{xushanhuang@kaist.ac.kr}

\author[Kang]{Moon-Jin Kang}
\address[Moon-Jin Kang]{\newline Department of Mathematical Sciences \newline Korea Advanced Institute of Science and Technology, Daejeon  34141, Republic of Korea}
\email{moonjinkang@kaist.ac.kr}

\author[Kim]{Jeongho Kim}
\address[Jeongho Kim]{\newline Department of Applied Mathematics, \newline Kyung Hee University, 1732 Deogyeong-daero, Giheung-gu, Yongin-si, Gyeonggi-do 17104, Republic of Korea}
\email{jeonghokim@khu.ac.kr}

\author[Lee]{Hobin Lee}
\address[Hobin Lee]{\newline Department of Mathematical Sciences \newline Korea Advanced Institute of Science and Technology, Daejeon  34141, Republic of Korea}
\email{lcuh11@kaist.ac.kr}

\begin{document}
	\newtheorem{theorem}{Theorem}[section]
	\newtheorem{lemma}{Lemma}[section]
	\newtheorem{corollary}{Corollary}[section]
	\newtheorem{proposition}{Proposition}[section]
	\newtheorem{remark}{Remark}[section]
	\newtheorem{definition}{Definition}[section]
	
	\renewcommand{\theequation}{\thesection.\arabic{equation}}
	\renewcommand{\thetheorem}{\thesection.\arabic{theorem}}
	\renewcommand{\thelemma}{\thesection.\arabic{lemma}}
	\newcommand{\bbr}{\mathbb R}
	\newcommand{\bbz}{\mathbb Z}
	\newcommand{\bbn}{\mathbb N}
	\newcommand{\bbs}{\mathbb S}
	\newcommand{\bbp}{\mathbb P}
	\newcommand{\bbt}{\mathbb T}
	\newcommand{\<}{\langle}
	\renewcommand{\>}{\rangle}
	\newcommand{\e}{\varepsilon}
	\newcommand{\pa}{\partial}
	\newcommand{\tU}{\widetilde{U}}
	\newcommand{\tu}{\widetilde{u}}
	\newcommand{\tv}{\widetilde{v}}
	\newcommand{\tw}{\widetilde{w}}

	\newcommand{\rd}{\partial}
	\newcommand{\na}{\nabla}
	
	\newcommand{\infi}{\infty}
	\newcommand{\R}{\mathbb{R}}	
	\newcommand{\sX}{\dot{\bold{X}}}
	\newcommand{\bB}{\bold{B}}
	\newcommand{\bG}{\bold{G}}
	\newcommand{\bS}{\bold{S}}
	\newcommand{\bD}{\bold{D}}
	\newcommand{\bX}{\bold{X}}
	\newcommand{\bY}{\bold{Y}}
	\newcommand{\cB}{\mathcal{B}}
	\newcommand{\cG}{\mathcal{G}}
	\newcommand{\cS}{\mathcal{S}}
	\newcommand{\cD}{\mathcal{D}}
	
	\newcommand{\tvR}{\tilde{v}^R}
	\newcommand{\tvS}{\tilde{v}^S}
	\newcommand{\tvX}{\tilde{v}^{\bold{X}}}
	\newcommand{\tvRX}{(\tilde{v}^R)^{\bold{X}}}
	
	\newcommand{\tuR}{\tilde{u}^R}
	\newcommand{\tuS}{\tilde{u}^S}
	\newcommand{\tuX}{\tilde{u}^{\bold{X}}}
	
	\newcommand{\pv}{p(v)}
	\newcommand{\tuRX}{(\tilde{v}^R)^{\bold{X}}}
	\newcommand{\tp}{\tilde{p}}
	\newcommand{\tpv}{p(\tilde{v))}}
	\newcommand{\tpvX}{p(\tilde{v}^{\bold{X}})}
	\newcommand{\tUS}{\tilde{U}^S}
	\newcommand{\UX}{U^{\bold{X}}}
	\newcommand{\tUX}{\tU^{\bold{X}}}
	\newcommand{\norm}[1]{\left\lVert#1\right\rVert}
	\newcommand{\beq}{\begin{equation}}
\newcommand{\eeq}{\end{equation}}

\newcommand{\eps}{\varepsilon }
	
	
	
	
	\thanks{\textbf{Acknowledgment.} X. Huang was partially supported by the National Research Foundation of Korea  (NRF-2019R1A5A1028324). M.-J. Kang and H. Lee were partially supported by the National Research Foundation of Korea  (NRF-2019R1C1C1009355  and NRF-2019R1A5A1028324). }

	\begin{abstract} 
		We consider the compressible barotropic Navier-Stokes equations in a half-line and study the time-asymptotic behavior toward the outgoing viscous shock wave. Precisely, we consider the two boundary problems:  impermeable wall and inflow problems, where the velocity at the boundary is given as a constant state. For both problems, when the asymptotic profile determined by the prescribed constant states at the boundary and far-fields is a viscous shock,   
		we show that the solution asymptotically converges to the shifted viscous shock profiles uniformly in space, under the condition that initial perturbation is small enough in $H^1$ norm. 		Since our method works on the physical variables, we do not require that the anti-derivative variables belong to $L^2$ space as in \cite{HMS03,MM99}. Moreover, for the inflow case, we remove the assumption $\gamma\le 3$ in \cite{HMS03}.
		Our results are based on the method of $a$-contraction with shifts, as the first extension of the method to the boundary value problems.
	\end{abstract}
	
	\maketitle
	
	\tableofcontents
	
	\section{Introduction}
	We consider the one-dimensional compressible barotropic Navier-Stokes(NS) equations in a half-line $\bbr_+=(0,\infty)$:
	\begin{align}
	\begin{aligned}\label{eq:NS-eulerian}
	&{\rho}_t+({\rho}{u})_{{x}}=0,\quad t\in\R_+,\quad {x}\in\R_+,\\
	&({\rho}{u})_t+({\rho}{u}^2+{p})_{{x}}=\mu{u}_{{x}{x}},
	\end{aligned}
	\end{align}
	subject to the boundary condition
	\begin{equation}\label{eq:boundary}
	({\rho},{u})|_{{x}=0}=(\rho_-,u_-), \quad t>0,
	\end{equation}
	and the far-field condition
	\begin{equation}\label{eq:initial}
	({\rho},{u}) \to(\rho_+,u_+),\quad\mbox{as}\quad {x}\to\infty, \quad t>0.
	\end{equation}
	The initial data $({\rho},{u})|_{{t}=0}=(\rho_0,u_0)$ satisfies the two conditions \eqref{eq:boundary}-\eqref{eq:initial}.
	Here, ${\rho}$ and ${u}$ denote the density and velocity of the fluid, and ${p}={\rho}^\gamma$ is the pressure with adiabatic constant $\gamma> 1$. For simplicity, $\mu$ is the viscosity coefficient, which will be chosen as $\mu=1$.
	
	The initial-boundary value problem \eqref{eq:NS-eulerian}--\eqref{eq:initial} is classified into three classes. When $u_->0$, which implies that the fluid flows through the boundary ${x}=0$, it is called the {\it inflow problem}; when $u_-<0$, that is the flow recede through the boundary, the problem is called the {\it outflow problem}; finally, when $u_-=0$, there is no flow through the boundary and therefore, the problem is called {\it impermeable wall problem}. In the case of the impermeable wall problem, the condition ${\rho}(t,0)=\rho_-$ should be removed and the boundary condition \eqref{eq:boundary} becomes simply ${u}(t,0)=u_-$. In the present paper, we consider the impermeable wall and inflow problems.\\
	
	\noindent $\bullet$ (Impermeable wall problem):
	We first focus on the impermeable wall problem. In this case, the NS equations \eqref{eq:NS-eulerian} together with the boundary condition ${u}_-=0$ can be represented in terms of the Lagrangian mass coordinates $(t,x)$, as in \cite{MM99} : (Here, we still keep the variable $x$ to denote the mass variable) 
	
	\begin{equation}\label{eq:NS}
	\begin{aligned}
		&v_t-u_x=0,\qquad t\in \mathbb{R}_+, \quad x\in \mathbb{R}_+,\\
		&u_t+p(v)_x=\left(\frac{u_x}{v}\right)_x, 
	\end{aligned}
	\end{equation}
	subject to the boundary and far-field conditions:
	\begin{equation}\label{boundary_impermeable}
	u(t,0)=u_-=0,\quad (v,u)\to(v_+,u_+)\quad\mbox{as}\quad x\to\infty, \quad t>0,
	\end{equation} 
	which is also satisfied by initial data $(v,{u})|_{{t}=0}=(v_0,u_0)$.
	Here, $v=1/\rho>0$ denotes the specific volume, $v_+:=1/\rho_+$ and the pressure is now given as $p(v)=v^{-\gamma}$. \\
	
	We are interested in the large-time behavior of the solution $(v,u)$ to the NS equations \eqref{eq:NS}. In particular, we consider the case when $u_+<0$, in which the solution is expected to converge toward the outgoing viscous shock (or equivalently, 2-shock) \cite{MM99}. The viscous shock profile uniquely exists up to shift if the values $(v_{\pm},u_{\pm})$ satisfy the Rankine-Hugoniot(RH) conditions:
	\begin{align}
	\begin{aligned}\label{RH}
	&-\sigma (v_+-v_-)-(u_+-u_-)=0,\\
	&-\sigma (u_+-u_-)+(p(v_+)-p(v_-))=0.
	\end{aligned}
	\end{align}
	In the case of impermeable wall condition $u_-=0$, for any given $v_+>0$ and $u_+<0$, the values of $v_-$ and the shock speed $\sigma$ are determined from the RH conditions \eqref{RH} as
	\begin{equation}\label{RH_impermeable}
	v_+^{-\gamma+1}\left(1-\left(\frac{v_-}{v_+}\right)^{-\gamma}\right)\left(1-\frac{v_-}{v_+}\right)=-u_+^2,\quad \sigma= \frac{-u_+}{v_+-v_-},
	\end{equation}
	with $v_-<v_+$. For these far field states $(v_\pm,u_\pm)$, the viscous shock wave $(\tv(\zeta),\tu(\zeta))$ with $\zeta=x-\sigma t$ is given by the solution to the following ODEs:
	\begin{equation}\label{eq:NS-shock}
	\begin{aligned}
		&-\sigma(\widetilde{v})'-(\widetilde{u})'=0,\\
		&-\sigma(\widetilde{u})'+p(\widetilde{v})'=\left(\frac{(\widetilde{u})'}{\widetilde{v}}\right)', \\
		&(\widetilde{v}, \widetilde{u})(-\infty)=(v_-, u_-), \quad (\widetilde{v}, \widetilde{u})(+\infty)=(v_+, u_+).
	\end{aligned}
	\end{equation}
	The first goal of the present paper is to attain the convergence of the solution $(v(t,x),u(t,x))$ to the impermeable wall problem \eqref{eq:NS}-\eqref{boundary_impermeable} towards the shifted viscous shock profile satisfying \eqref{eq:NS-shock}.\\
	
	\noindent $\bullet$ (Inflow problem): On the other hand, for the case of the inflow problem, the NS system  \eqref{eq:NS-eulerian} is transformed to the Lagrangian mass coordinates, as in \cite{HMS03,M01} : (We still keep the variable $x$ to denote the mass variable)
	
	\begin{align}
	\begin{aligned} \label{eq:NS_inflow}
	&v_t -u_x = 0,\quad t>0,\quad x>\sigma_-t,\\
	&u_t+p(v)_x=\left(\frac{u_x}{v}\right)_{x},
	\end{aligned}
	\end{align}
	subject to the boundary and far-field conditions:
	\beq\label{bdy-in}
	(v,u)(t,x)|_{x=\sigma_-t}=(v_-,u_-),\quad (v,u)\to(v_+,u_+)\quad \mbox{as}\quad x\to\infty,\quad t>0,
	\eeq
	where $\sigma_-=-\frac{u_-}{v_-}<0$ and $v_\pm:=1/\rho_\pm$. We fix boundary $x=\sigma_-t$ of the system \eqref{eq:NS_inflow} by using a change of variable $\xi=x-\sigma_-t>0$, and so \eqref{eq:NS_inflow} is rewritten as
	\begin{align}
	\begin{aligned}\label{eq:NS_inflow_2}
	&v_t -\sigma_- v_\xi -u_\xi = 0,\quad t>0,\quad \xi>0,\\
	&u_t-\sigma_-u_\xi +p(v)_\xi=\left(\frac{u_\xi}{v}\right)_{\xi},
	\end{aligned}
	\end{align}
	with the boundary and far-field conditions:
	\begin{equation}\label{boundary_inflow}
	(v,u)(t,0)=(v_-,u_-),\quad (v,u)\to(v_+,u_+)\quad \mbox{as}\quad \xi\to\infty,\quad t>0,
	\end{equation}
	which is also satisfied by initial data $(v,{u})|_{{t}=0}=(v_0,u_0)$.

	 Different from the impermeable wall problem, the boundary value $v_-$ should be imposed and the values $(v_\pm,u_\pm)$ satisfy the same RH conditions \eqref{RH} and the entropy condition $v_-<v_+$ and $u_->u_+$, where now the shock speed $\sigma$ becomes
	\begin{equation}\label{shock_speed_inflow}
	\sigma = \sqrt{-\frac{p(v_+)-p(v_-)}{v_+-v_-}}>0.
	\end{equation}
  	We are also interested in the time-asymptotic behavior of solutions to the inflow problem \eqref{eq:NS_inflow_2}-\eqref{boundary_inflow} (and \eqref{eq:NS_inflow}-\eqref{bdy-in}) toward the outgoing viscous shock $(\tv(\zeta),\tu(\zeta))$ with $\zeta=x-\sigma t=\xi-(\sigma-\sigma_-)t$ satisfying \eqref{eq:NS-shock}. 
	For that, as in \cite{HMS03}, we consider when the boundary state $(v_-,u_-)$ belongs to the subsonic region $\Omega_{sub}:= \{ (v,u) \in \bbr_+\times\bbr_+ ~|~ |u|< v \sqrt{-p'(v)}  \}$ and the far-field state  $(v_+,u_+)$ lies on the 2-shock curve starting from $(v_-,u_-)$, that is, the RH condition \eqref{RH} and the entropy condition $v_-<v_+$, $u_->u_+$ hold.
	For this situation, it is expected from a heuristic argument \cite{HMS03,M01} that solutions $(v(t,\xi),u(t,\xi))$ to \eqref{eq:NS_inflow_2} asymptotically converge to the viscous shock profile $(\tu(\xi-(\sigma-\sigma_-)t),\tv(\xi-(\sigma-\sigma_-)t))$ up to shift.	For the rest of the paper, we interpret the solution $(v,u)$ as a function of $(t,x)$ for the impermeable wall problem, and as a function of $(t,\xi)$ for the inflow problem.\\
	
	\subsection{Literature review}
	Large-time behaviors of the one-dimensional compressible NS equations have been extensively studied after the pioneering work of Matsumura-Nishihara \cite{MN85,MN86}, Goodman \cite{G86}, Liu \cite{L85}, and Szepessy-Xin \cite{SZ93}. Although there are plenty of results when the spatial domain is the whole line $\bbr$, we focus on reviewing the previous results on the large-time behavior of the compressible NS equations for the case of half-line $\bbr_+$. For the impermeable wall problem, the asymptotic behavior of the solution can be classified into two cases, namely the outgoing viscous shock and rarefaction wave, which are investigated in \cite{MM99} and \cite{MN99} respectively. However, the patterns for asymptotic behaviors of the inflow and outflow problems are much more complicated, due to the presence of the boundary layer solution, and they are classified depending on the states at the boundary and far-field \cite{M01}. In \cite{MN01}, the asymptotic convergence towards both the rarefaction wave and the superposition of the boundary layer solution and the rarefaction wave are established. On the other hand, the cases of viscous shock wave and its superposition with the boundary layer solution are considered in \cite{HMS03}. Finally, for the case of outflow problems, the asymptotic stability of the boundary layer solution and the superposition of rarefaction wave and the boundary layer solution were studied in \cite{HQ09, KNZ03, KZ08}. We also refer to \cite{HLS10, QW09,QW11} and references therein for the large-time behaviors of the Navier-Stokes-Fourier equations on the half-line.
	
	Note that the results in \cite{HMS03,MM99} regarding the stability of the viscous shock wave are based on the anti-derivative method which requires that the anti-derivatives of the physical variables belong to $L^2$. 
	Recently, the stability results on the viscous shock wave and its composition with the other elementary waves are investigated in \cite{HKK23,KVW23,KVW-NSF} by using the method of $a$-contraction with shifts developed in \cite{KV16}, which does not introduce the anti-derivative variables. Instead, it directly uses the original perturbation itself, and so does not need additional assumption on the initial data of the anti-derivative variable.\\	

	\subsection{Main results}
	
	We now state the main results on the global existence and  large-time behavior of the impermeable wall problem and the inflow problem.
	
	\begin{theorem}[Impermeable wall problem] \label{thm:impermeable}
		Assume $\gamma>1$.  For a given constant state $(v_+,u_+)$ with $v_+>0$ and $u_+<0$, there exists $\delta_0,\e_0,\beta>0$ such that the following holds. \\
		For any $(v_-,u_-=0)$ satisfying \eqref{RH_impermeable} such that $|v_+-v_-|\sim |u_+-u_-|=|u_+|<\delta_0$, let $(\tv(x-\sigma t),\tu(x-\sigma t))$ be the viscous 2-shock wave satisfying \eqref{eq:NS-shock} where the shock speed $\sigma$ is given in \eqref{RH_impermeable}. Let $(v_0,u_0)$ be any initial data such that
		\begin{equation}\label{initial_condi}
		\|(v_0,u_0)-(v_+,u_+)\|_{L^2(\beta,\infty)}+\|(v_0,u_0)-(v_-,u_-)\|_{L^2(0,\beta)}+\|(\pa_xv_0,\pa_x u_0)\|_{L^2(\R_+)}<\e_0.
		\end{equation}
		Then, the impermeable wall problem \eqref{eq:NS}--\eqref{boundary_impermeable} subject to the initial data $(v_0,u_0)$ admits a unique global-in-time solution $(v,u)(t,x)$ as follows: there exists a Lipschitz shift  $t\mapsto \bX(t)$ such that
		\begin{align*}
		&(v,u)(t,x)-(\tv,\tu)(x-\sigma t -\bX(t)-\beta)\in C(0,\infty;H^1(\R_+)),\\
		&u_{xx}(t,x)-\tu_{xx}(x-\sigma t-\bX(t)-\beta)\in L^2(0,\infty;L^2(\R_+)).
		\end{align*}
		Moreover, 
		\begin{equation}\label{long_time_impermeable}
		\lim_{t\to\infty}\sup_{x\in\R_+}\left|(v,u)(t,x)-(\tv,\tu)(x-\sigma t-\bX(t)-\beta)\right|=0
		\end{equation}
		and $\lim_{t\to\infty}|\dot{\bX}(t)|=0$.
	\end{theorem}

For the inflow problem, a similar statement holds. The only difference is that we use the variable $\xi=x-\sigma_-t$ instead of the variable $x$.	
	
	\begin{theorem}[Inflow problem]\label{thm:inflow}
		Assume $\gamma>1$. For a given constant state $(v_-,u_-)\in\Omega_{sub}:= \{ (v,u) \in \bbr_+\times\bbr_+ ~|~ |u|< v \sqrt{-p'(v)}  \}$, there exists $\delta_0,\e_0,\beta>0$ such that the following holds. 
		
		For any $(v_+,u_+)\in\bbr_+\times\bbr$ satisfying \eqref{RH} and the entropy condition $v_-<v_+$ and $u_->u_+$ such that $|v_+-v_-|\sim |u_+-u_-|<\delta_0$, let $(\tv(\xi-(\sigma-\sigma_-) t),\tu(\xi-(\sigma-\sigma_-) t))$ be the viscous 2-shock wave satisfying \eqref{eq:NS-shock} where the shock speed $\sigma$ is given in \eqref{shock_speed_inflow}. Let $(v_0,u_0)$ be any initial data such that
		\begin{equation}\label{infini}
		\|(v_0,u_0)-(v_+,u_+)\|_{L^2(\beta,\infty)}+\|(v_0,u_0)-(v_-,u_-)\|_{L^2(0,\beta)}+\|(\pa_xv_0,\pa_x u_0)\|_{L^2(\R_+)}<\e_0.
		\end{equation}
		Then, the inflow problem \eqref{eq:NS_inflow_2}--\eqref{boundary_inflow} subject to the initial data $(v_0,u_0)$ admits a unique global-in-time solution $(v,u)(t,\xi)$ as follows: there exists a Lipschitz shift  $t\mapsto \bX(t)$ such that
		\begin{align*}
		&(v,u)(t,\xi)-(\tv,\tu)(\xi-(\sigma-\sigma_-) t -\bX(t)-\beta)\in C(0,\infty;H^1(\R_+)),\\
		&u_{\xi\xi} (t,\xi) -\tu_{\xi\xi}(\xi-(\sigma-\sigma_-) t-\bX(t)-\beta)\in L^2(0,\infty;L^2(\R_+)).
		\end{align*}
		Moreover, 
		\begin{equation}\label{long_time_inflow}
		\lim_{t\to\infty}\sup_{\xi\in\R_+}\left|(v,u)(t,\xi)-(\tv,\tu)(\xi-(\sigma-\sigma_-) t-\bX(t)-\beta)\right|=0
		\end{equation}
		and $\lim_{t\to\infty}|\dot{\bX}(t)|=0$.
	\end{theorem}

	\begin{remark} 
	 As in \cite{HMS03,MM99}, we consider the time-asymptotic stability of viscous shock under small perturbation in $H^1$ norm. For that, we consider the situation where the initial data $(v_0,u_0)$ are given in a neighborhood of $(\tv,\tu)(x-\beta)$ for some large constant $\beta>0$ as in \eqref{initial_condi} and \eqref{infini}. Then, our results imply the long-time behavior toward the asymptotic profile $(\tv,\tu)(x-\sigma t-\beta)$ up to a dynamical shift $\bX(t)$ where the asymptotic proflie is far away from the boundary.\\
	 Indeed, in the above results, the decay estimate $\lim_{t\to\infty}|\dot{\bX}(t)|=0$ implies
		$$
		\lim_{t\rightarrow+\infty}\frac{\bX(t)}{t}=0,
		$$
		which means that the shift function $\bX(t)$ grows at most sub-linearly as $t\to\infty$. Thus, the dynamically shifted wave $(\tv,\tu)(x-\sigma t-\bX(t)-\beta)$ time-asymptotically tends to the original wave $(\tv,\tu)(x-\sigma t-\beta)$ shifted by the constant $\beta$.
	\end{remark}
	
	\begin{remark} 
	 Theorems \ref{thm:impermeable} and \ref{thm:inflow} do not require that the initial anti-derivatives of the perturbations belong to $L^2$, which is imposed by the previous results \cite{HMS03,MM99}. In particular, for the inflow problem, Theorem \ref{thm:inflow} holds for any $\gamma>1$, which removes the assumption $\gamma\le 3$ of \cite{HMS03}.
	Our results handle a shock wave of small jump strength, while the results of \cite{HMS03,MM99} are for a shock wave of moderate strength. However, we do not need the smallness assumption on $u_-$ that is crucially used in \cite{HMS03} for the inflow problem.
	\end{remark}

	\subsection{Main ideas for the proof}
	As we mentioned above, the main tool for proving the results in Theorem \ref{thm:impermeable} and Theorem \ref{thm:inflow} is the {\it method of $a$-contraction with shift}, which was developed in \cite{KV16} to study the stability of the shock wave. Consider a system with an entropy $\eta$. For example, the NS system \eqref{eq:NS} is equipped with a natural entropy (as a mechanical energy of the system) $\eta(U)=\frac{u^2}{2}+Q(v)$, where
	\[U=(v,u),\quad Q(v) = \frac{v^{1-\gamma}}{\gamma-1}.\]
	Then, one of the natural quantities to measure the difference between the state $U$ and the reference state $\overline{U}$ is the {\it relative entropy} defined by
	\[\eta(U|\overline{U})=\eta(U)-\eta(\overline{U})-\nabla\eta(\overline{U})(U-\overline{U}),\]
	which is always nonnegative due to the convexity of $\eta$ with respect to $U$ and it vanishes if and only if $U=\overline{U}$. 
In \cite{KV21,KV-Inven}, the method of $a$-contraction with shift was used to show the contraction of large perturbations of a viscous shock. More precisely, if $U$ is any large NS solution, and $\overline{U}$ is a viscous shock, then
\[
\frac{d}{dt}\int_\Omega a^{\bold{X}}  \eta(U|\overline{U}^{\bold{X}})\,dx\le 0,
\]	
where some weight function $a$ and the viscous shock are shifted by some dynamical shift $\bold{X}(t)$.
	This contraction property plays a crucial role in the inviscid limit problem as in \cite{KV-Inven,KV-2shock}.
The property also successfully derives the large-time behavior of the composite waves of NS \cite{HKK23,KVW23} or Naiver-Stokes-Fourier systems \cite{KVW-NSF}. In the same context, the goals of the present paper are to carry out an establishment of the method of $a$-contraction with shift to the boundary value problems:
	the impermeable wall and inflow problems. \\
	
	However, when applying the method to the boundary problems, we would naturally encounter the following main difficulties.\\
	First, when applying the Poincar\'e-type inequality of Lemma \ref{Poincare} to the main terms in \eqref{terms_YJP} that are localized by derivative of weight and viscous shock, we use the change of variable $x\mapsto y$ as in \eqref{ydef} that is defined the viscous shock. However, the range $[y_0,1]$ of the new variable $y$ depends on time, contrary to the whole space case $\bbr$ (as in \cite{KV21,KVW23}), since  the viscous shock connects the left state $(v_-,u_-)$ at $x=-\infty$ to the right state $(v_+,u_+)$ at $x=+\infty$ whereas the solution stays only on the half space $\bbr_+$. 
	To overcome it, we use the smallness of $|\bX(t)|$ to get $0<y_0\le  Ce^{-C\delta\beta}$ as in \eqref{bbtu}, by which $y_0 \approx 0$ for $\beta$ large enough.  
	 At this point, we also use the important fact that the optimal constant $\frac{1}{2}$ of the Poincar\'e-type inequality does not depend on the range of $y$ as in Lemma \ref{Poincare}.\\
	Second, we should control the size of boundary values arising in using integration by parts, as in Lemmas \ref{lem:boundary}, \ref{lem:v-H1} and \ref{lem:boundary_inflow}. For that, we use the interpolation inequality, the a prior bounds and some cancellation structure.

	\subsection{Organization of the paper}
	The paper is organized as follows. In Section \ref{sec:2}, we provide preliminaries on the estimates for the relative quantities, viscous shock wave, and the Poincar\'e-type inequality, which are used in the later analysis. Section \ref{sec:3} presents a statement of the a priori estimates on $H^1$-perturbation of the solution. Based on them, we also provide the proof of the main theorems. Then, we provide detailed proof of the a priori estimate on the impermeable wall problem in Section \ref{sec:4} and Section \ref{sec:5}. Finally, we provide the proof of the a priori estimate on the inflow problem in Section \ref{sec:6}, highlighting the difference in the proof compared to the impermeable wall problem.

	\section{Preliminaries}\label{sec:2}
	\setcounter{equation}{0}
	
	We first present several estimates on the relative quantities and viscous shock wave, which will be used in the later analysis. We also provide the Poincar\'e-type inequality for general domain.
	
	\subsection{Estimates on the relative quantities}
	For any function $F:(0,\infty)\to \bbr$ and $v,\overline{v}\in(0,\infty)$, we define the relative quantity $F(v|\overline{v})$ of $F$ as
	\[F(v|\overline{v}):=F(v)-F(\overline{v})-F'(\overline{v})(v-\overline{v}).\]
	In particular, we consider the relative quantities for the pressure $p(v)=v^{-\gamma}$ and the internal energy $Q(v)=\frac{v^{1-\gamma}}{\gamma-1}$, and present the lower and upper bound estimates on them. Since the proof of the following lemma can be found in \cite{KV21}, we omit the proof for simplicity.
	\begin{lemma}\label{lem:relative}
		Let $\gamma>1$ and $v_+$ be given constants. Then, there exist constants $C$ and $\delta_*$ such that the following assertions hold:
		\begin{enumerate}
			\item For any $v,\overline{v}$ satisfying $0<\overline{v}<2v_+$ and $0<v<3v_+$,
			\[|v-\overline{v}|^2\le CQ(v|\overline{v}),\quad |v-\overline{v}|^2\le Cp(v|\overline{v}).\]
			\item For any $v,\overline{v}>\frac{v_+}{2}$,
			\[|p(v)-p(\overline{v})|\le C|v-\overline{v}|.\]
			\item For any $0<\delta<\delta_*$ and any $v,\overline{v}>0$ satisfying $|p(v)-p(\overline{v})|<\delta$ and $|p(\overline{v})-p(v_+)|<\delta$,
			\begin{align*}
			&p(v|\overline{v})\le\left(\frac{\gamma+1}{2\gamma p(\overline{v})}+C\delta\right)|p(v)-p(\overline{v})|^2,\\
			&Q(v|\overline{v})\ge \frac{p(\overline{v})^{-\frac{1}{\gamma}-1}}{2\gamma}|p(v)-p(\overline{v})|^2-\frac{1+\gamma}{3\gamma^2}p(\overline{v})^{-\frac{1}{\gamma}-2}(p(v)-p(\overline{v}))^3,\\
			&Q(v|\overline{v})\le\left(\frac{p(\overline{v})^{-\frac{1}{\gamma}-1}}{2\gamma}+C\delta\right)|p(v)-p(\overline{v})|^2.
			\end{align*}
		\end{enumerate}
	\end{lemma}

	\subsection{Viscous shock wave}
	Recall that the profile of the viscous shock wave satisfies the ODEs \eqref{eq:NS-shock}, and the existence, uniqueness, and properties of the viscous shock wave are now well understood. The following lemma summarizes the main properties of the viscous shock wave, which will be used in the later analysis. We refer to \cite{G86,MN85} for the proof of the lemma.
	
	\begin{lemma}\label{lem:viscous_shock}
		For a given right-end state $(v_+,u_+)$, there exists a constant $C>0$ such that the following holds. For any left-end state $(v_-,u_-)$ connected with $(v_+,u_+)$ via 2-shock curve, there exists a unique solution $(\tv(\zeta),\tu(\zeta))$ to \eqref{eq:NS-shock}. Let $\delta$ be the strength of the shock defined as $\delta:=|u_+-u_-|\sim|v_+-v_-|$. Then, we have
		\[\tv_\zeta>0,\quad \tu_\zeta<0,\]
		and
		\begin{align}
		\begin{aligned}\label{shock_prop}
		&|\tv(\zeta)-v_{\pm}|, |\tu(\zeta)-u_{\pm}|\le C\delta e^{-C\delta|\zeta|},\quad \pm \zeta>0,\\
		&|(\tv_\zeta,\tu_\zeta)|\le C\delta^2 e^{-C\delta|\zeta|},\quad \zeta\in\R,\\
		&|(\tv_{\zeta\zeta},\tu_{\zeta\zeta})|\le C\delta|(\tv_\zeta,\tu_\zeta)|,\quad \zeta\in\R.
		\end{aligned}
		\end{align}
	\end{lemma} 

	\subsection{Poincar\'e-type inequality}
	One of the main tools for attaining asymptotic stability result is the Poincar\'e-type inequality. 
	The optimal constant $\frac{1}{2}$ of the  Poincar\'e-type inequality in \cite{KV21} is independent of the size of domain as follows. This is useful in our analysis.

	\begin{lemma}\label{Poincare}
		For any $c<d$ and function $f : \left[c, d\right]  \longrightarrow  \mathbb{R}$ satisfying $\int_{c}^d (y-c)(d-y)|f'(y)|^2dy<\infty$,   
		\[\int_{c}^d \left|f(y)-\frac{1}{d-c}\int_{c}^d f(y)dy\right|^2dy\le \frac{1}{2}\int_{c}^d(y-c)(d-y)|f'(y)|^2dy.\]
	\end{lemma}
	\begin{proof}
	First, we recall the Poincar\'e-type inequality in \cite{KV21}. For any $g:[0,1]\to\R$ satisfying $\int_0^1 z(1-z)|g'(z)|^2\,dz<\infty$, we have
	\begin{equation}\label{poincare-original}
		\int_0^1 \left|g-\int_0^1 g\,dz\right|^2\,dz\le \frac{1}{2}\int_0^1 z(1-z)|g'|^2\,dz.
	\end{equation}
	Then, we use change of variables $y=c +(d-c)z$ to observe that
	\begin{align*}
		\frac{1}{d-c}\int_c^d f(y)\,dy=\int_0^1g(z)\,dz,\quad g(z):=f(y)=f(c+(d-c)z).
	\end{align*}
	Therefore we use \eqref{poincare-original} to obtain
	\begin{align*}
		\int_c^d&\left|f(y)-\frac{1}{d-c}\int_c^d f(y)\,dy\right|^2\,dy\\
		&=(d-c)\int_0^1 \left|g(z)-\int_0^1g(z)\,dz\right|^2\,dz\le \frac{d-c}{2}\int_0^1 z(1-z)|g'(z)|^2\,dz\\
		&=\frac{d-c}{2}\int_c^d \frac{y-c}{d-c}\frac{d-y}{d-c}(d-c)^2|f'(y)|^2\frac{dy}{d-c}\\
		&=\frac{1}{2}\int_c^d(y-c)(d-y)|f'(y)|^2\,dy.
	\end{align*}
	\end{proof}

	\section{A priori estimate and proof of the main theorem}\label{sec:3}
	\setcounter{equation}{0}
	In this section, we state a priori estimates for the $H^1$-perturbation between the solution and the viscous shock wave, and based on them, we prove the large-time behavior of the impermeable and inflow problem towards the viscous shock waves.
	
	\subsection{Local existence of solutions}
	We first ensure that the impermeable wall problem admits a unique local-in-time solution.
	
	\begin{proposition}\label{prop:local}
		For any constant $\beta>0$, let $\underline{v}$ and $\underline{u}$ be smooth monotone functions such that
		\[(\underline{v}(x),\underline{u}(x))=(v_+,u_+),\quad\mbox{for}\quad x\ge\beta,\quad \underline{v}(0)>0.\]
		For any constants $M_0$, $M_1$, $\underline{\kappa}_0$, $\overline{\kappa}_0$, $\underline{\kappa}_1$, and $\overline{\kappa}_1$ with $0<M_0<M_1$ and $0<\underline{\kappa}_1<\underline{\kappa}_0<\overline{\kappa}_0<\overline{\kappa}_1$, there exists a constant $T_0>0$ such that if
		\begin{align*}
			&\|(v_0-\underline{v},u_0-\underline{u})\|_{H^1(\R_+)}\le M_0,\\
			&0<\underline{\kappa}_0\le v_0(x)\le \overline{\kappa}_0,\quad x\in\R_+,
		\end{align*}
		the impermeable wall problem \eqref{eq:NS}--\eqref{boundary_impermeable} has a unique solution $(v,u)$ on $[0,T_0]$ such that
		\begin{align*}
			v-\underline{v}\in C([0,T_0];H^1(\R_+)),\quad u-\underline{u}\in C([0,T_0];H^1(\R_+))\cap L^2(0,T_0;H^2(\R_+)),
		\end{align*}
		and
		\[\|(v-\underline{v},u-\underline{u})\|_{L^\infty(0,T_0;H^1(\R_+))}\le M_1.\]
		Moreover,
		\[\underline{\kappa}_1\le v(t,x)\le \overline{\kappa}_1,\quad u(t,0)=0,\quad \forall t>0,\quad x\in\R_+.\]
	\end{proposition}

	\begin{proof}
		Since the local existence can be proved using the standard iterative method by considering the following sequence of functions $\{(v^n,u^n)\}_{n=1}^\infty$:
		\begin{align*}
		&v^{n+1}_t-u^{n+1}_x=0,\quad t>0,\quad x>0,\\
		&u^{n+1}_t+p(v^n)_x=\left(\frac{u^{n+1}_x}{v^n}\right)_x,\quad u^{n+1}(t,0)=0,
		\end{align*}
		we omit the proof.
	\end{proof}

	Similarly, we also obtain the local existence for the inflow problem, that is presented by the same statement as in Proposition \ref{prop:local} by replacing the boundary condition $u(t,0)=0$ by $(v,u)(t,0)=(v_-,u_-)$.

	\subsection{Construction of shift}
	As we mentioned in the introduction, it is expected that the viscous shock wave should be shifted to attain the stability estimate. In this part, we explicitly construct the shift function. For the impermeable wall problem, we first define the weight function $a=a(\zeta)$ as
	\begin{equation}\label{def_a}
	a(\zeta) := 1-\frac{\tu(\zeta)}{\sqrt{\delta}},\quad \zeta = x-\sigma t
	\end{equation}
	where $\delta = |u_+-u_-|=|u_+|$ denotes the strength of shock. We note that the weight function $a$ satisfies $1\le a\le 1+\sqrt{\delta}<\frac{3}{2}$ for small enough $\delta$, and 
	\begin{equation}\label{a_derivative}
	a'(\zeta)=-\frac{ \tu'(\zeta)}{\sqrt{\delta}}=\frac{\sigma \tv'(\zeta)}{\sqrt{\delta}}>0,\quad \mbox{and}\quad |a'(\zeta)|\sim \frac{\tv'(\zeta)}{\sqrt{\delta}}.
	\end{equation} 
	Then, we define the shift function $\bX:\bbr_+\to\bbr$ as a solution to the following ODE:
	\begin{equation}
	\begin{aligned}\label{def_shift}
		\sX(t)=-\frac{M}{\delta}\Bigg(&\int_{\bbr_+}a^{\bX,\beta}(\zeta)\tu^{\bX,\beta}_\zeta(\zeta)(u(t,x)-\tu^{\bX,\beta}(\zeta))\,dx\\
		\quad &+\frac{1}{\sigma}\int_{\bbr_+}a^{\bX,\beta}(\zeta)\pa_x[p(\tv^{\bX,\beta}(\zeta))](v(t,x)-\tv^{\bX,\beta}(\zeta))\,dx\Bigg),\quad \bX(0)=0,
	\end{aligned}
	\end{equation}
	where $M, \beta>0$ are positive constants which will be chosen later. Here, for any function $f:\bbr\to\bbr$, we use the abbreviated notation
	\[f^{\bX,\beta}(\cdot):=f(\cdot-\bX(t)-\beta).\]

	Similarly, the weight function and shift function for the inflow problem are defined exactly in the same manner, except that we use the variable $\xi=x-\sigma_-t$ instead of $x$. More precisely, define the weight function $a$ for the inflow problem as
	\[a(\zeta):=1+\frac{u_--\tu(\zeta)}{\sqrt{\delta}},\quad \zeta = \xi-(\sigma-\sigma_-)t,\]
	where $\delta=|u_+-u_-|$ and the same estimates in \eqref{a_derivative} hold. We define the shift function $\bX$ for the inflow problem as
	\begin{equation}
	\begin{aligned}\label{def_shift_inflow}
		\sX(t)=-\frac{M}{\delta}\Bigg(&\int_{\R_+}a^{\bX,\beta}(\zeta)\tu_\zeta^{\bX,\beta}(\zeta)(u(t,\xi)-\tu^{\bX,\beta}(\zeta))\,d\xi\\
		&+\frac{1}{\sigma}\int_{\R_+}a^{\bX,\beta}(\zeta)\pa_\xi[ p(\tv^{\bX,\beta}(\zeta))](v(t,\xi)-\tv^{\bX,\beta}(\zeta))\,d\xi\Bigg),\quad \bX(0)=0.
	\end{aligned}
	\end{equation}
		As in \cite[Lemma 3.3]{KVW23}, we ensure the existence of Lipschitz solution to the above ODEs \eqref{def_shift} and \eqref{def_shift_inflow} under the the condition that $v$ has positive upper and lower bounds, and $u$ is bounded, which is guaranteed by Proposition \ref{prop:local} and the a priori assumption \eqref{apriori_small}. \\
		
	At this moment, it is unclear why the shifts are defined in the above manner. However, this becomes clear in Section \ref{sec:4}, where we apply the method of $a$-contraction with shift.
	
	\subsection{A priori estimates}
	
	Now, we are ready to present the main proposition on a priori estimates for both the impermeable wall problem and the inflow problem, which are the key estimates for obtaining the large-time behaviors of the systems. In the following proposition, for the   impermeable problem, the shifted viscous shock wave $(\tv^{\bX,\beta},\tu^{\bX,\beta})$ is understood as a function of $(t,x)$ :
	\[\tv^{\bX,\beta}(t,x)=\tv(x-\sigma t-\bX(t)-\beta),\quad \tu^{\bX,\beta}(t,x)=\tu(x-\sigma t-\bX(t)-\beta).\]
For 	the inflow problem, the shifted viscous shock wave is understood as a function of $(t,\xi)$ :
\[\tv^{\bX,\beta}(t,\xi)=\tv(\xi-(\sigma-\sigma_-) t-\bX(t)-\beta),\quad \tu^{\bX,\beta}(t,\xi)=\tu(\xi-(\sigma-\sigma_-) t-\bX(t)-\beta).\]
	
	\begin{proposition}\label{prop:main}
		For a given $(v_+,u_+)\in\bbr_+\times\bbr$, there exist positive constants $C_0,\delta,\e$, and $\beta$ such that the following holds.
		Suppose that $(v,u)$ is the solution to the impermeable wall problem \eqref{eq:NS}--\eqref{boundary_impermeable} (resp.  inflow problem \eqref{eq:NS_inflow_2}--\eqref{boundary_inflow})  on $[0,T]$ for some $T>0$, and $(\tv,\tu)$ is the viscous shock profile satisfying \eqref{eq:NS-shock} with $\delta=|u_+|$ (resp. $\delta=|u_+-u_-|$), and the shift $\bX$ is defined in \eqref{def_shift} (resp. \eqref{def_shift_inflow}). Suppose that
		\begin{align*}
			&v-\tv^{\bX,\beta}\in C([0,T];H^1(\R_+)),\\
			&u-\tu^{\bX,\beta}\in C([0,T];H^1(\R_+))\cap L^2(0,T;H^2(\R_+)),
		\end{align*}
		and
		\begin{equation}\label{apriori_small}
			\|v-\tv^{\bX,\beta}\|_{L^\infty(0,T;H^1(\R_+))}+\|u-\tu^{\bX,\beta}\|_{L^\infty(0,T;H^1(\R_+))}\le \e.
		\end{equation}
		Then, for all $0\le t\le T$,
		\begin{align}
			\begin{aligned}\label{apriori_impermeable}
				\sup_{t\in[0,T]}&\left[\|v-\tv^{\bX,\beta}\|_{H^1(\R_+)}+\|u-\tu^{\bX,\beta}\|_{H^1(\R_+)}\right]\\
				&\quad +\sqrt{\int_0^t (\delta|\sX|^2+G_1+ G^S+D_{v_1}+D_{u_1}+D_{u_2})\,ds}\\
				&\le C_0\left(\|v_0-\tv^{\bX,\beta}(0,\cdot)\|_{H^1(\R_+)}+\|u_0-\tu^{\bX,\beta}(0,\cdot)\|_{H^1(\R_+)}\right)+C_0e^{-C\delta\beta},
			\end{aligned}
		\end{align}
		where $C_0$ is independent of $T$ and
		\begin{align}
			\begin{aligned}\label{good_terms}
				&G_1:=\int_{\R_+}a_x\left|p(v)-p(\tv^{\bX,\beta})-\frac{u-\tu^{\bX,\beta}}{2C_*}\right|^2\,dx,\quad G^S:=\int_{\R_+}|\tu^{\bX,\beta}_x||u-\tu^{\bX,\beta}|^2\,dx,\\
				&D_{v_1}:=\int_{\R_+}|(p(v)-p(\tv^{\bX,\beta}))_{x}|^2\,dx,\\
				&D_{u_1}:=\int_{\R_+}|(u-\tu^{\bX,\beta})_x|^2\,dx,\quad D_{u_2}:=\int_{\R_+}|(u-\tu^{\bX,\beta})_{xx}|^2\,dx,
			\end{aligned}
		\end{align}
where $x$ is replaced by $\xi$ for the inflow problem.\\		
 In particular, for all $0\le t \le T$,
 \begin{equation} \label{bddx12}
|\dot{\bX}(t)|\le C_0\lVert (v-\widetilde{v}^{\bX,\beta})(t,\cdot) \rVert_{L^\infty(\mathbb{R}_+)}.
\end{equation}
	\end{proposition}

	\subsection{Proof of Theorem \ref{thm:impermeable} and Theorem \ref{thm:inflow}}
Based on Propositions \ref{prop:local} and \ref{prop:main}, we use the continuation argument to prove  the global-in-time existence of perturbations. We also use Proposition \ref{prop:main} to prove the long-time behavior. Those proofs are similar to that of the previous articles (e.g. \cite{HKK23,KVW23}). Therefore, we present them in Appendix, and complete the proofs of Theorem \ref{thm:impermeable} and Theorem \ref{thm:inflow}.\\
	
	Therefore, it only remains to prove the main Proposition \ref{prop:main}. 
	We will present detailed proofs of Proposition \ref{prop:main} for the impermeable case in Section \ref{sec:4} and Section \ref{sec:5}. Since the proof of Proposition \ref{prop:main} for the inflow case shares a lot of parts with the impermeable case, we will only provide estimates for some new terms in Section \ref{sec:6}.\\
	

	\section{Relative Entropy estimates for the impermeable wall problem}\label{sec:4}
	\setcounter{equation}{0}
	
	In this section, we obtain estimates on the $L^2$-norms of the perturbation for the impermeable wall problem, by using the method of $a$-contraction with shift.  
	In what follows, we suppress the dependency on the shift $\bX$ and $\beta$, that is, we will use the following concise notation with no confusion:
	\begin{align*}
	&a(t,x)=a(x-\sigma t-\bX(t)-\beta),\\
	&\tv(t,x) = \tv^{\bX,\beta}(t,x)=\tv(x-\sigma t-\bX(t)-\beta),\\
	&\tu(t,x) = \tu^{\bX,\beta}(t,x)=\tu(x-\sigma t-\bX(t)-\beta).
	\end{align*}
In the remaining part, $C$ denotes a positive $O(1)$-constant that may change from line to line, but is independent of the parameters $\delta, \eps, \beta$ and the time $T$.	
	
\noindent The goal of this section is to prove the following lemma.
	
	\begin{lemma}\label{lem:main}
		Under the hypothesis of Proposition \ref{prop:main}, there exists a positive constant $C$ such that for all $t\in[0,T]$,
		\begin{align}
		\begin{aligned} \label{mainlemmaend}
		&\|v-\tv\|_{L^2(\R_+)}^2+\|u-\tu\|_{L^2(\R_+)}^2+\int_0^t (\delta|\sX(s)|^2+G_1+G^S+D_{u_1})\,ds\\
		&\quad \le C\left(\|v_0-\tv(0,\cdot)\|_{L^2(\R_+)}^2+\|u_0-\tu(0,\cdot)\|^2_{L^2(\R_+)}\right)+Ce^{-C\delta \beta}+C\varepsilon^2\int_{0}^t\|(u-\widetilde{u})_{xx}\|_{L^2(\mathbb{R}_+)}^{2}ds,
		\end{aligned}
		\end{align}
		where $G_1$, $G^S$, and $D_{u_1}$ are the terms defined in \eqref{good_terms}.
	\end{lemma}

	\subsection{Weighted relative entropy estimate}
	To prove Lemma \ref{lem:main}, we use the celebrated relative entropy method introduced by Dafermos and DiPerna \cite{D96,D79}. To this end, we rewrite the system \eqref{eq:NS} into the general form of viscous hyperbolic conservation laws:
	\begin{equation} \label{eq:general_hyperbolic} 
		\rd_t U+\rd_x A(U)=\rd_x \left( M(U) \rd_x \nabla \eta (U) \right),
	\end{equation}
	where
	\[U:=\begin{pmatrix}
	v \\
	u
	\end{pmatrix}, \quad
	A(U):=\begin{pmatrix}
	-u\\
	p(v)
	\end{pmatrix}.\]
	Here the entropy $\eta$ of the system \eqref{eq:general_hyperbolic} is given by $\eta(U): ={u^2 \over 2}+Q(v)$,  where $Q(v)={v^{-\gamma +1} \over \gamma+1}$,  i.e.,  $Q'(v)=-p(v)$ and the matrix $M(U)$ is given by
	\[M(U):=\begin{pmatrix}
	0 &0\\
	0 &{1 \over v}
	\end{pmatrix}.\]
	Similarly, the shifted viscous shock wave
	\[\tU(t,x):=\begin{pmatrix}
	\tv(t,x)\\\tu(t,x)
	\end{pmatrix}\]
	satisfies a similar system given as
	\begin{equation}\label{eq:general_shock}
		\rd_t \tU+\rd_x A(\tU)=\rd_x \left( M(\tU) \rd_x \nabla \eta (\tU) \right)-\dot{\bX}\rd_x\tU.
	\end{equation}
	To estimate the difference between $U$ and $\tU$, we use the relative entropy as a measurement of the difference. The relative entropy functional is defined as
	\[
		\eta(U | V)=\eta(U)-\eta(V)-\nabla \eta(V)(U-V),
	\]
	and the relative flux is given by
	\[
		A(U|V)=A(U)-A(V)-\nabla A(V)(U-V).
	\]
	Finally, let $G(\cdot;\cdot)$ be the relative entropy flux defined as
	\[
		G(U; V)=G(U)-G(V)-\nabla \eta(V) (A(U)-A(V)),
	\]
	where $G$ is the entropy flux of $\eta$,  i.e.,  $\rd_i G(U)=\sum_{k=1}^2 D_k \eta(U) \rd_i A_k (U)$ for $i=1,2$. For the system \eqref{eq:general_hyperbolic}, a straight computation yields
	\begin{equation} \label{relative_quantities}
		\begin{split}
			&\eta(U|\tU)={|u-\tu|^2 \over 2}+Q(v|\tv), \quad A(U|\tU)=\begin{pmatrix}
				0\\
				p(v|\tv)
			\end{pmatrix},\quad G(U;\tU)=(p(v)-p(\tv))(u-\tu).
		\end{split}
	\end{equation}
	
	Below, we will estimate the relative entropy, weighted by the function $a$ in \eqref{def_a}, between the solution $U(t,x)$ of \eqref{eq:general_hyperbolic} and the shifted viscous shock wave $\tU(t,x)$ of \eqref{eq:general_shock}.
	
	\begin{lemma}\label{lem:rel_ent_est}
	Let $a$ be the weight function defined in \eqref{def_a}, $U$ be a solution to \eqref{eq:general_hyperbolic} and $\tU$ be the shifted shock wave satisfying \eqref{eq:general_shock}. Then,
		\begin{equation}\label{est-1}
			{d \over dt}\int_{\mathbb{R}_+} a(t,x)\eta (U(t,x))|\tU(t,x)) d x=\sX(t)\bY+\mathcal{J}^{\textup{bad}}-\mathcal{J}^{\textup{good}}+\mathcal{P},
		\end{equation}
		where
		\begin{equation}\label{terms_YJP}
			\begin{split}
				\bY&:= -\int_{\mathbb{R}_+} a_x \eta (U|\tU)\,dx+\int_{\mathbb{R}_+} a\nabla^2 \eta (\tU)(\tU)_x(U-\tU)\,dx , \\
				\mathcal{J}^{\textup{bad}}&:=\int_{\mathbb{R}_+} a_x(p(v)-p(\tv))(u-\tu)\,dx-\int_{\mathbb{R}_+} a \widetilde{u}_xp(v | \tv) \,dx\\
				&\quad-\int_{\mathbb{R}_+} a_x {u-\tu \over v} \rd_x (u-\tu)\,dx+ \int_{\mathbb{R}_+} a_x(u-\tu)(v-\tv){\rd_x \tu \over v \tv}\,dx\\
				&\quad+\int_{\mathbb{R}_+} a\rd_x (u-\tu) {v-\tv \over v \tv} \rd_x \tu  \,dx,\\
				\mathcal{J}^{\textup{good}}&:={\sigma \over 2}\int_{\mathbb{R}_+} a_x|u-\tu|^2\,dx+\sigma \int_{\R_+} a_x Q(v|\tv)\,dx\\
				&\quad+ \int_{\mathbb{R}_+} {a \over v} |\rd_x (u-\tu)|^2\,dx,\\
				\mathcal{P}&:=\left[-a\widetilde{u}(p(v)-p(\widetilde{v}))+a\widetilde{u}\frac{(u_x-\widetilde{u}_x)}{v}-a\widetilde{u}(v-\widetilde{v})\frac{\widetilde{u}_x}{v\widetilde{v}}\right]_{x=0}.
			\end{split}
		\end{equation}
	\end{lemma}
	\begin{remark} Since $a_x>0$, the terms in $\mathcal{J}^{\textup{good}}$ have positive signs as good terms, while the signs of the terms in $\mathcal{J}^{\textup{bad}}$ are indefinite. Moreover, the terms in $\mathcal{P}$ come from the boundary value, when we take the integration-by-parts from $x=0$ to $x=+\infty$.
	\end{remark}
	\begin{proof}

		Following  the same computations as in \cite[Lemma 2.3]{KV21}, it holds from \eqref{eq:general_hyperbolic} and \eqref{eq:general_shock} that
		\[{d \over dt}\int_{\mathbb{R}_+} a \eta (U|\tU) dx=\sX(t) \bY -\sigma\int_{\mathbb{R}_+}a_x\eta(U|\widetilde{U})dx+\sum_{i=1}^5I_i,\]
		where
		\begin{equation}\label{I_i}
			\begin{split}
				&I_1:=-\int_{\mathbb{R}_+} a \rd_x G(U; \tU) dx ,\\
				&I_2:=-\int_{\mathbb{R}_+} a \rd_x \nabla \eta (\tU) A(U|\tU) d x ,\\
				&I_3:=\int_{\mathbb{R}_+} a \left( \nabla \eta (U)-\nabla \eta(\tU)\right) \rd_x \left(M(U) \rd_x(\nabla \eta(U)-\nabla \eta(\tU) \right) dx,\\
				&I_4:= \int_{\mathbb{R}_+} a \left( \nabla \eta(U)-\nabla \eta (\tU)\right) \rd_x \left( (M(U)-M(\tU))\rd_x \nabla \eta (\tU)\right) d x,\\
				&I_5:=\int_{\mathbb{R}_+} a(\nabla \eta )(U|\tU)\rd_x\left( M(\tU)\rd_x \nabla \eta(\tU)\right)dx.
			\end{split}
		\end{equation}
		Although the overall estimates for $I_i$ for $i=1,2,\ldots, 5$ are the same as in the previous literature, e.g., \cite{KV21}, we present the estimates on them for the readers' convenience, especially taking care of the boundary terms.\\
		
		\noindent $\bullet$ (Estimate of $I_1$): We use the definition of $G(\cdot;\cdot)$ in \eqref{relative_quantities}, boundary condition $u(t,0)=0$, and integration-by-parts to derive
		\begin{align*}
			I_1&=\int_{\mathbb{R}_+} a_x G(U;\tU) d x+\left[aG(U;\tU)\right]_{x=0}\\
			&=\int_{\mathbb{R}_+} a_x(p(v)-p(\tv))(u-\tu)dx-a(t,0)\widetilde{u}(t,0)(p(v(t,0))-p(\widetilde{v}(t,0))).
		\end{align*}
		\noindent $\bullet$ (Estimate of $I_2$): Using the definition of $A(U|\tU)$ in \eqref{relative_quantities}, we have
		\begin{align*}
			I_2&=-\int_{\mathbb{R}_+} a \tu_x p(v|\tv) dx.
		\end{align*}
		\noindent $\bullet$ (Estimates of $I_3$ and $I_4$): As in the estimate of $I_1$, we use integration-by-parts to obtain
		\begin{align*}
			I_3&=\int_{\mathbb{R}_+} a (u-\tu) \rd_x \left( {1 \over v} \rd_x (u-\tu) \right) d x\\
			&=-\int_{\mathbb{R}_+} {a \over v} |\rd_x (u-\tu)|^2 dx-\int_{\R_+} a_x {u-\tu \over v} \rd_x (u-\tu) dx\\
			& \quad +a(t,0)\widetilde{u}(t,0)\frac{(u_x(t,0)-\widetilde{u}_x(t,0))}{v(t,0)},
		\end{align*}
		and
		\begin{align*}	
			I_4&=\int_{\mathbb{R}_+} a(u-\tu) \rd_x \left( {\tv-v \over v \tv} \rd_x \tu \right) dx\\
			&=\int_{\mathbb{R}_+} a_x(u-\tu)(v-\tv){\rd_x \tu \over v \tv} \,dx+\int_{\mathbb{R}_+} a \rd_x (u-\tu) {v-\tv \over v \tv} \rd_x \tu \,dx\\
			&\quad -a(t,0)\widetilde{u}(t,0)(v(t,0)-\widetilde{v}(t,0))\frac{\widetilde{u}_x(t,0)}{v(t,0)\widetilde{v}(t,0)}.
		\end{align*}
		\noindent$\bullet$ (Estimate of $I_5$): Finally, it can be directly obtained that $I_5=0$.\\
		
		\noindent Therefore, combining the estimates on $I_i$ for $i=1,2,\ldots,5$, we obtain
		\begin{align*}
			&{d \over dt} \int_{\mathbb{R}_+} a (t,x) \eta(U(t,x)|\tU (t,x))\,dx\\
			&=\sX(t) \bY+\int_{\mathbb{R}_+} a_x(p(v)-p(\tv))(u-\tu)dx-\int_{\mathbb{R}_+} a \tu_x p(v|\tv)\,dx\\
			&\quad -\int_{\mathbb{R}_+} a_x {u-\tu \over v} \rd_x (u-\tu) dx + \int_{\mathbb{R}_+} a_x(u-\tu)(v-\tv){\rd_x \tu \over v \tv} dx+\int_{\mathbb{R}_+} a \rd_x (u-\tu) {v-\tv \over v \tv} \rd_x \tu \,dx\\
			&\quad -{\sigma \over 2}\int_{\mathbb{R}_+} a_x|u-\tu|^2\,dx-\sigma \int_{\mathbb{R}_+} a_x Q(v|\tv)\,dx-\int_{\mathbb{R}_+} {a \over v} |\rd_x (u-\tu)|^2\,dx\\
			&\quad+\left[-a\widetilde{u}(p(v)-p(\widetilde{v}))+a\widetilde{u}\frac{(u_x-\widetilde{u}_x)}{v}-a\widetilde{u}(v-\widetilde{v})\frac{\widetilde{u}_x}{v\widetilde{v}}\right]_{x=0}\\
			&=\sX(t)\bY+\mathcal{J}^{\textup{bad}}-\mathcal{J}^{\textup{good}}+\mathcal{P}.
		\end{align*}
	\end{proof} 
	
	\subsection{Maximization on $p(v)-p(\tv)$}
	
	Our strategy is to control the bad terms $\mathcal{J}^{\textup{bad}}$ on the right-hand side of \eqref{est-1} by using the good terms $\mathcal{J}^{\textup{good}}$. However, the most troubling one among the terms in $\mathcal{J}^{\textup{bad}}$ is
	\[\int_{\R_+}a_x(p(v)-p(\tv))(u-\tu)\,dx,\]
	where the perturbation $p(v)-p(\tv)$ and $u-\tu$ are coupled. Therefore, to decouple those perturbations, we will rewrite $\mathcal{J}^{\textup{bad}}$ into the maximized representation in terms of $p(v)-p(\tv)$. We will use the following lemma, which is exactly the same as \cite[Lemma 4.3]{HKKL_pre}.
	\begin{lemma}\cite[Lemma 4.3]{HKKL_pre}\label{lem:quad}
		 For any $\delta>0$ small enough, let $C_*$ be the constant as
		\[C_*:=\frac{1}{2}\left(\frac{1}{\sigma_l}-\sqrt{\delta}\frac{\gamma+1}{\gamma}\frac{1}{p(v_-)}\right), \quad \sigma_l:=\sqrt{-p'(v_-)}.
		\]
		Then we have
		\begin{align*}
			-\int_{\mathbb{R}_+}& a \tu_x p(v|\tv) \, dx-\sigma \int_{\mathbb{R}_+} a_x Q(v|\tv) \, dx\\
			&\le -C_*\int_{\mathbb{R}_+}a_x|p(v)-p(\tv)|^2\,dx\\
			&\quad +C\delta \int_{\mathbb{R}_+} a_x\big|p(v)-p(\tv)\big|^2 \, d x+C\int_{\mathbb{R}_+} a_x\big|p(v)-p(\tv)\big|^3 \, d x.
		\end{align*}
	\end{lemma}
	Using Lemma \ref{lem:quad} and the quadratic structure with respect to $p(v)-p(\tv)$, we can represent $\mathcal{J}^{\textup{bad}}-\mathcal{J}^{\textup{good}}$ in the another form.
	\begin{lemma}\label{lem:badgood}
		For $\mathcal{J}^{\textup{bad}}$ and $\mathcal{J}^{\textup{good}}$ defined in \eqref{terms_YJP}, we have
		\[\mathcal{J}^{\textup{bad}}-\mathcal{J}^{\textup{good}} \le \mathcal{B}-\mathcal{G},\]
		where
		\begin{align*}
			\mathcal{B}&:=\frac{1}{4 C_*}\int_{\mathbb{R}_+} a_x |u-\widetilde{u}|^2 \, dx+C\delta \int_{\mathbb{R}_+} a_x\big|p(v)-p(\tv)\big|^2 \, d x+C\int_{\mathbb{R}_+} a_x\big|p(v)-p(\tv)\big|^3 \, d x \\
			&\quad-\int_{\mathbb{R}_+} a_x\frac{u-\tu}{v} \partial_x (u-\tu) dx+\int_{\mathbb{R}_+}a_x\frac{(u-\tu)(v-\tv) \partial_x \tu}{v \tv} dx+\int_{\mathbb{R}_+} a  \partial_x (u-\tu) \frac{(v-\tv)\partial_x \tu}{v \tv} \,dx\\
			\mathcal{G}&:=C_* \int_{\mathbb{R}_+} a_x \left| p(v)-p(\tv)-\frac{u-\tu}{2C_*}\right|^2 \, d x+ \frac{\sigma}{2} \int_{\mathbb{R}_+} a_x |u-\tu|^2 \, dx+\int_{\mathbb{R}_+} \frac{a}{v} |\partial_x (u-\tu)|^2 \, dx.
		\end{align*} 
	Again, since $\sigma>0$, $a_x >0$, and $a>0$,  the terms in $\mathcal{G}$ have positive sign.
	\end{lemma}

	\begin{proof}
		Let $J_1$ and $J_2$ be the first and second term of $\mathcal{J}^{\textup{bad}}$, and let $J_3$ be the second term of $\mathcal{J}^{\textup{good}}$:
		\begin{align*}
			J_1:=\int_{\mathbb{R}_+} a_x(p(v)-p(\tv))(u-\tu)\,dx,\quad J_2:=-\int_{\mathbb{R}_+} a \widetilde{u}_xp(v | \tv)\,dx,\quad J_3:=\sigma \int_{\mathbb{R}_+} a_x Q(v|\tv) \,dx.
		\end{align*}
		Then, we use Lemma \ref{lem:quad} to derive
		\begin{align}
			\begin{aligned}\label{J123}
				J_1+J_2-J_3
				&\le \int_{\mathbb{R}_+} a_x (p(v)-p(\widetilde{v}))(u-\widetilde{u})\,dx-C_*\int_{\mathbb{R}_+}a_x|p(v)-p(\tv)|^2\,dx\\
				&\quad +C\delta \int_{\mathbb{R}_+} a_x\big|p(v)-p(\tv)\big|^2 \, d x+C\int_{\mathbb{R}_+} a_x\big|p(v)-p(\tv)\big|^3 \, d x.
			\end{aligned}
		\end{align}
		However, since the first two terms in the right-hand side of \eqref{J123} can be rewritten by using the quadratic structure of $p(v)-p(\tv)$ as
		\begin{align*}
			&\int_{\mathbb{R}_+} a_x (p(v)-p(\widetilde{v}))(u-\widetilde{u})dx-C_* \int_{\mathbb{R}_+} a_x |p(v)-p(\widetilde{v})|^2 dx\\
			&=\int_{\mathbb{R}_+} a_x \left[ -C_* \left( (p(v)-\widetilde{p}(v))^2-\dfrac{(p(v)-p(\widetilde{v}))(u-\widetilde{u})}{C_*} + \dfrac{(u-\widetilde{u})^2}{4 C^*} \right) + \dfrac{(u-\widetilde{u})^2}{4C_*}  \right] d x\\
			&=\int_{\mathbb{R}_+} a_x \left[ -C_* \left( (p(v)-p(\widetilde{v}))- \dfrac{u-\widetilde{u}}{2 C_*} \right)^2 + \dfrac{(u-\widetilde{u})^2}{4C_*}  \right] d x.
		\end{align*}
		Since the other terms in $\mathcal{J}^{\textup{bad}}$ and $\mathcal{J}^{\textup{good}}$ are unchanged, we obtain the desired inequality.
	\end{proof}

\subsection{Bound for shift}
Applying Lemma \ref{lem:viscous_shock} to \eqref{bddx12}, and using $|(\widetilde{v})'|\sim|(\widetilde{u})'|$, we have
\begin{align*}
		|\dot{\bX}(t)| &\le \frac{C}{\delta}  \left(\lVert (v-\tv(t,\cdot) \rVert_{L^\infty(\mathbb{R}_+)}+\lVert (u-\tu(t,\cdot) \rVert_{L^\infty(\mathbb{R}_+)}\right) \int_{\mathbb{R}_+} |\tv_{x} | dx \\
		&\le C \left(\lVert (v-\tv(t,\cdot) \rVert_{L^\infty(\mathbb{R}_+)}+\lVert (u-\tu(t,\cdot) \rVert_{L^\infty(\mathbb{R}_+)}\right),
		\end{align*}   
		which yields \eqref{bddx12}.
	\subsection{Estimate on the boundary terms}
	Compared to the previous results \cite{HKK23,KVW23} on the large-time behaviors on $\R$, the most prominent difference is the presence of the boundary terms $\mathcal{P}$. Since the terms in $\mathcal{P}$ are pointwise values at $x=0$, they should be treated separately. In the following lemma, we provide the time integration of the boundary terms can be controlled by the constant shift $\beta$ and the second-order term of $u-\tu$.
	
	\begin{lemma} \label{lem:boundary}
		There exists $C>0$ independent of $\delta$ such that
		\[
		\left|\int_{0}^t\mathcal{P}\, d s\right|  \le Ce^{-C\delta \beta}+C\varepsilon^2\int_{0}^t\|(u-\widetilde{u})_{xx}\|_{L^2(\mathbb{R}+)}^{2}\,ds 
		\]
		for all $t\in \left[0, T\right]$.
	\end{lemma}
	\begin{proof}
First, using \eqref{bddx12}, the a priori assumption \eqref{apriori_small} with Sobolev inequality and the smallness of $\eps$, we have
\[
|\dot{\bX}(t)|\le C\eps \le \frac{\sigma}{2}, \quad t\le T,
\]
and so
\[
|{\bX}(t)| \le \frac{\sigma}{2} t, \quad t\le T,
\]
which yields
 \begin{equation} \label{bddx0}
-\sigma t -\bX (t) -\beta \le - \frac{\sigma}{2} t -\beta <0, \quad t\le T.
\end{equation}
            We split the boundary terms into three parts as
		\begin{align*}
		\mathcal{P}&=\left[-a\widetilde{u}(p(v)-p(\widetilde{v}))\right]_{x=0}+\left[a\widetilde{u}\frac{(u_x-\widetilde{u}_x)}{v}\right]_{x=0}+\left[-a\widetilde{u}(v-\widetilde{v})\frac{\widetilde{u}_x}{v\widetilde{v}}\right]_{x=0} \\
		&=:\mathcal{P}_1+\mathcal{P}_2+\mathcal{P}_3.
		\end{align*}     
		We first use \eqref{shock_prop} and \eqref{bddx0} to have 
            \beq \label{bbtu}
		|\widetilde{u}(t,0)|= |\tu(-\sigma t-\bX(t)-\beta)-u_-| \le C\delta e^{-C\delta|-\sigma t-\bX(t)-\beta|}=C\delta e^{-C\delta t}e^{-C\delta \beta}.
		\eeq
		Then we use \eqref{apriori_small} to estimate $\mathcal{P}_1$ as 
		\begin{align*}
		\left|\int_{0}^{t}\mathcal{P}_1\,ds\right|\le C\lVert p(v)-p(\widetilde{v}) \rVert_{L^\infty(\mathbb{R}_+)} \int_{0}^{t} |\widetilde{u}(s,0)|\,ds \le C\varepsilon e^{-C\delta \beta}.    
		\end{align*}
		Using the interpolation inequality, Young's inequality and \eqref{apriori_small}, we obtain
		\begin{align*}
		\left|\int_{0}^{t}\mathcal{P}_2(U)\,ds\right| &= \left|\int_{0}^{t}\frac{a}{v}(u-\widetilde{u})(u-\widetilde{u})_{x}\big|_{x=0} ds \right| \\
		& \leq C\int_{0}^{t} |\widetilde{u}(s,0)| \|(u-\widetilde{u})_{x}\|_{L^\infty(\mathbb{R}_+)}ds \\
		& \leq C\int_{0}^{t} |\widetilde{u}(s,0)|^{\frac{4}{3}} ds + C\int_{0}^{t} \|(u-\widetilde{u})_{x}\|_{L^2(\mathbb{R}_+)}^{2}\|(u-\widetilde{u})_{xx}\|_{L^2(\mathbb{R}_+)}^{2}ds \\
		& \leq C\delta^{\frac{1}{3}}e^{-C\delta\beta} + C\varepsilon^{2}\int_{0}^{t}\|(u-\widetilde{u})_{xx}\|_{L^2(\mathbb{R}_+)}^{2}ds.
	\end{align*} 
Likewise, since
\[
|\widetilde{u}_x(t,0)|= |\tu'(-\sigma t-\bX(t)-\beta)| \le C\delta^2  e^{-C\delta t}e^{-C\delta \beta},
\]
we have
		\begin{align*}
		\left|\int_{0}^{t}\mathcal{P}_{3}\,ds\right|\le C\|v-\tv\|_{L^\infty}\int_0^t\left| \Big[\frac{\tu\tu_x}{v\tv}\Big]_{x=0}\right|\,ds\le C\varepsilon \delta^2e^{-C\delta\beta}.
		\end{align*}		
		Combining the estimates on $\mathcal{P}_i$ for $i=1,2,3$, we finally derive
		\[\left|\int_0^t \mathcal{P}\,ds\right|\le Ce^{-C\delta\beta}+C\varepsilon^2\int_{0}^t\|(u-\widetilde{u})_{xx}\|_{L^2(\mathbb{R}_+)}^{2}\,ds.\]
	\end{proof}
	
	\subsection{Proof of Lemma \ref{lem:main}}\label{sec:proof_main_1}
	We now provide the proof of Lemma \ref{lem:main}. It follows from the results of Lemma \ref{lem:rel_ent_est} and Lemma \ref{lem:badgood} that 
	\begin{equation}\label{est-rel_ent-2}
		{d \over dt}\int_{\mathbb{R}_+} a \eta (U | \tU) \,dx \le \sX(t)\bY+\mathcal{B}-\mathcal{G}+\mathcal{P}.
	\end{equation}
	We split the bad terms $\mathcal{B}$ and the good terms $\mathcal{G}$ as follows:
	\[\mathcal{B}:=\sum_{i=1}^6 B_i,\quad \mathcal{G}:=G_1+G_2+D,\]
	where
	\begin{align*}
		&B_1:=\frac{1}{4 C^*} \int_{\mathbb{R}_+} a_x  |u-\tu|^2\,dx,\quad B_2:=-\int_{\mathbb{R}_+} a_x {u-\tu \over v} \rd_x (u-\tu)\,dx,\\
		&B_3:= \int_{\mathbb{R}_+} a_x(u-\tu)(v-\tv){\rd_x \tu \over v \tv} \,dx,\quad B_4:=\int_{\mathbb{R}_+} a \rd_x (u-\tu) {v-\tv \over v \tv} \rd_x \tu \,dx,\\
		&B_5:=C\delta \int_{\mathbb{R}_+} a_x\left|p(v)-p(\tv)\right|^2 \, d x,\quad B_6:=C\int_{\mathbb{R}_+} a_x\left|p(v)-p(\tv)\right|^3 \, d x, 
	\end{align*}
	and
	\begin{align*}
		&G_1:= C_* \int_{\mathbb{R}_+} a_x \left|p(v)-p(\tv)-\frac{u-\tu}{2C^*} \right|^2\,d x,\quad G_2:={\sigma \over 2} \int_{\mathbb{R}_+} a_x |u-\tu|^2\,dx,\\
		&D:= \int_{\mathbb{R}_+} {a \over v} |\rd_x (u-\tu)|^2\,dx.
	\end{align*}
	Moreover, we also split the term $\bY$ as
	\begin{align*}
		\bY&=-\int_{\mathbb{R}_+} a_x \eta (U | \tU)\, dx+\int_{\mathbb{R}_+} a \nabla^2 \eta (\tU)\tU_x (U-\tU)\,dx\\
		&=-\int_{\mathbb{R}_+} a_x \left( {|u-\tu|^2 \over 2}+Q(v|\tv)\right)\,dx+\int_{\mathbb{R}_+} a \tu_x (u-\tu)\,dx-\int_{\mathbb{R}_+} a p'(\tv) \tv_x ( v-\tv)\,dx\\
		&=\sum_{i=1}^6 \bY_i,
	\end{align*}
	where
	\begin{align*}
		\bold{Y}_{1} &:=\int_{\mathbb{R}_+} a \tu_x (u-\tu) \,dx,\quad \bold{Y}_{2}:=\frac{1}{\sigma}\int_{\mathbb{R}_+} a p'(\tv)\tv_x (u-\widetilde{u}) \,dx,\\
		\bold{Y}_{3} &:=-\frac{1}{2}\int_{\mathbb{R}_+} a_x \Big( u-\tu -2C_* \big(p(v)-p(\tv) \big)\Big)\Big(u-\tu +2C_* \big(p(v)-p(\tv) \big)\Big) \, dx,\\
		\bold{Y}_{4} &:=- 2C_*^2 \int_{\mathbb{R}_+} a_x  \big(p(v)-p(\tv) \big)^2 \,dx-\int_{\mathbb{R}_+} a_x Q(v|\tv)\,dx,\\
		\bold{Y}_{5} &:=- \int_{\mathbb{R}_+} a p'(\tv) \tv_x \left(v-\tv+\frac{2C_*}{\sigma}(p(v)-p(\tv)) \right) \, d x, \\
		\bold{Y}_{6} &:= \int_{\mathbb{R}_+} a p'(\tv) \tv_x \frac{2C_*}{\sigma}\left(p(v)-p(\tv)-\frac{u-\tu}{2C_*}\right) \, d x.
	\end{align*}
	Now, the ODE in \eqref{def_shift} can be written as
	\[
		\sX(t)=-\frac{M}{\delta}(\bold{Y}_1+\bold{Y}_2),
	\]
	which implies
	\[
		\sX(t)\bY=-\frac{\delta}{M}|\dot{\bX}(t)|^2+{\dot{\bX}}(t)\sum_{i=3}^6 \bY_i.
	\]
	Therefore, the right-hand side of \eqref{est-rel_ent-2} can be written as follows:
	\begin{align}
		\begin{aligned} \label{decomml}
			{d \over dt}\int_{\mathbb{R}_+} a \eta (U | \tU) dx&=\underbrace{-\frac{\delta}{2M}| {\dot{\bold{X}}}|^2+B_1-G_2-{3 \over 4}D}_{=:\mathcal{R}_1}\\
			&\quad \underbrace{-\frac{\delta}{2M}| {\dot{\bold{X}}}|^2+\dot{\bold{X}}\sum\limits_{i=3}^6\bold{Y}_i+\sum\limits_{i=2}^6B_i-G_1-\frac{D}{4}}_{:=\mathcal{R}_2}+\mathcal{P}.
		\end{aligned}    
	\end{align}
	Observe that, the definition \eqref{def_shift} of shift $\bX$ is motivated by the above decomposition. Among the bad terms $B_i$, the first term $B_1$ is the worst term to be controlled, as the other bad terms are comparably small, thanks to the smallness of the perturbation or shock strength. Therefore, estimating $\mathcal{R}_1$ would be the main issue, which uses a careful analysis and sharp Poincar\'e-type inequality in Lemma \ref{Poincare}. After obtaining the estimates for $\mathcal{R}_1$, the estimate on $\mathcal{R}_2$ can be obtained in a straightforward manner.
	
	\subsubsection{Estimate of $\mathcal{R}_1$}
	We first focus on obtaining the following estimate on $\mathcal{R}_1$.
	\begin{lemma} \label{lem:est_R1}
	For sufficiently large $\beta>0$, there exists a constant $C_1>0$ independent of $\delta$ such that 
		\begin{equation}\label{est_R1}
			\mathcal{R}_1=-\frac{\delta}{2M}| {\dot{\bold{X}}}|^2+B_1-G_2-{3 \over 4}D\leq -C_1\int_{\mathbb{R}_+} |\widetilde{u}_x||u-\widetilde{u}|^2 dx=:-C_1G^S.
		\end{equation}
	\end{lemma}
	\begin{proof}
		Since $\bX$ is bounded by \eqref{bddx12}, for a fixed $t\in [0,T]$, we consider a change of variable $x\mapsto y$ as
		\beq\label{ydef}
		y:=\frac{u_--\tu(x-\sigma t-\bX(t)-\beta)}{\delta}=-\frac{\tu(x-\sigma t-\bX(t)-\beta)}{\delta}.
		\eeq
		Indeed, the map $x\mapsto y=y(x)$ is one-to-one and increasing function satisfying
		\[\frac{dy}{dx} = -\frac{\tu_x(x-\sigma t-\bX(t)-\beta)}{\delta}>0,\quad \lim_{x\to0}y=y_0,\quad \lim_{x\to+\infty}y=1,\]
		where
		\[y_0 := \frac{-\tu(-\sigma t-\bX(t)-\beta)}{\delta}>0.\]
		As in \eqref{bbtu}, we have
		\[y_0\le Ce^{-C\delta |\sigma t+\bX(t)+\beta|}\le Ce^{-C\delta\beta},\]
		which implies that $y_0$ can be chosen arbitrary small by taking large $\beta$. In particular, we choose large enough $\beta$ so that $y_0<\frac{1}{6}$. Moreover, we represent the perturbation $u-\tu$ in terms of the variable $y$ as
		\[f(y):=(u(t,\cdot)-\tu(\cdot-\sigma t-\bX(t)-\beta))\circ y^{-1}.\]
		We also note that the weight function $a$ defined in \eqref{def_a} implies that $a(t,x)=1+\sqrt{\delta} y$ and therefore $\partial_xa=\sqrt{\delta} (dy / dx)$. 		Below, we will get a sharp estimate for each term of $\mathcal{R}_1$ separately.\\
		
		In the sharp estimates, we will use the following $O(1)$-constants:
	\begin{equation}\label{def_O1_const}
	\sigma_l:=\sqrt{-p'(v_-)}, \quad \alpha_l:={\gamma+1 \over 2 \gamma \sigma_l p(v_-)},
	\end{equation}
	which are indeed independent of the shock strength $\delta$. 
	The following approximations for $\sigma_l$ will be importantly used:
	\begin{align}
	\begin{aligned}\label{sigma_bound}
	&|\sigma-\sigma_l| \leq C \delta,\quad\|\sigma_l^2-p'(\tv(\cdot))\|_{L^\infty(\mathbb{R}_+)} \leq C \delta,\quad \left\| {1 \over \sigma^2_l}-{p(\tv(\cdot))^{-{1\over \gamma}-1} \over \gamma} \right\|_{L^\infty(\mathbb{R}_+)} \leq C \delta.
	\end{aligned}
	\end{align}
				
		\noindent $\bullet$ (Estimate of  $-\frac{\delta}{2M}|\sX|^2$): We first represent $\bY_1$ and $\bY_2$ by using the new variable $y$ as
		\begin{align*}
			&\bY_1=\int_{\R_+}a\tu_x(u-\tu)\,dx=-\delta \int_{y_0}^1 (1+\sqrt{\delta}y)f\,dy\\
			&\bY_2=-\frac{1}{\sigma^2} \int_{\mathbb{R}_+} a p'(\tv)\tu_x (u-\tu) \, d x=\frac{\delta}{\sigma^2} \int_{y_0}^1 (1+\sqrt{\delta}y)p'(\tv) f \, dy.
		\end{align*}
		Therefore, we have
		\begin{equation} \label{est_Y1}
			\left| \bY_1+\delta \int_{y_0}^1 f dy \right| \leq \delta^{3/2} \int_{y_0}^1 |f| dy,
		\end{equation}
		and by using \eqref{sigma_bound},
		\begin{equation}\label{est_Y2}
			\begin{split}
				\left|\bold{Y}_2+\delta \int_{y_0}^1 f \,dy \right| &\leq \left| \bold{Y}_2 -\frac{\delta}{\sigma^2} \int_{y_0}^1 p'(\tv)f \,dy\right|+\left|\frac{\delta}{\sigma^2}\int_{y_0}^1 p'(\tv)f\,dy+\delta \int_{y_0}^1 f\, dy\right| \\
				&\leq \frac{\delta^{3/2}}{\sigma^2}\int_{y_0}^1 y |p'(\tv)| |f|dy+\frac{\delta}{\sigma^2} \int_{y_0}^1 \left|\sigma^2+p'(\tv)\right||f|\,dy \\
				&\leq C (\delta^{3/2}+\delta^2) \int_{y_0}^1 |f|\,dy\le C\delta^{3/2}\int_{y_0}^1|f|\,dy.
			\end{split}
		\end{equation}
		Combining the estimates \eqref{est_Y1} and \eqref{est_Y2}, we obtain
		\begin{align*}
			\left| \sX-2M\int_{y_0}^1 f\,dy\right|=\left| \sum_{i=1}^2 \frac{M}{\delta} \left( \bY_i+\delta \int_{y_0}^1 f\, dy \right) \right| 
			&\leq C \sqrt{\delta} \int_{y_0}^1 |f|\,dy,
		\end{align*}
		which yields
		\begin{equation*}
			\left(  \left| 2M \int_{y_0}^1 f\,dy  \right| -|\sX| \right)^2 \leq C\delta\int_{y_0}^1 |f|^2\,dy.
		\end{equation*}
		The above inequality and a simple inequality $\frac{p^2}{2}-q^2 \leq (p-q)^2$ for all $p,q \in \mathbb{R}$ yield
		\begin{equation*}
			2M^2 \left( \int_{y_0}^1f\,dy \right)^2-|\sX|^2 \leq C\delta \int_{y_0}^1 |f|^2\,dy.
		\end{equation*}
		Thus, we obtain
		\begin{equation}\label{xdot}
			-\frac{\delta}{2M}|\sX|^2 \leq -M \delta  \left( \int_{y_0}^1 f\,dy \right)^2+C\delta^2 \int_{y_0}^1 |f|^2\,dy.
		\end{equation}
		
		\noindent $\bullet$ (Estimates of $B_1$ and  $G_2$): Recall that $B_1$ and $G_2$ are defined as
		\begin{align*}
			&B_1=\frac{1}{4C_*}\int_{\mathbb{R}_+} a_x |u-\tu|^2 \, d x, \quad G_2=\frac{\sigma}{2} \int_{\mathbb{R}_+} a_x |u-\tu|^2 \, d x.
		\end{align*}
		Therefore,
		\begin{align*}
			B_1-G_2&=\left(\frac{1}{4C_*}-\frac{\sigma}{2}\right)\int_{\mathbb{R}_+}a_x|u-\tu|^2\,dx=-\left(\frac{1}{4C_*}-\frac{\sigma}{2}\right)\frac{1}{\sqrt{\delta}}\int_{\mathbb{R}_+}\tu_x|u-\tu|^2\,dx\\
			&=\sqrt{\delta}\left(\frac{1}{4C_*}-\frac{\sigma}{2}\right)\int_{y_0}^1 |f|^2\,dy,
		\end{align*}
		where the constant $C_*$ defined in Lemma \ref{lem:quad} can be written as
		\begin{equation*}
			C_*= \frac{1}{2 \sigma_l} -\sqrt{\delta} \alpha_l \sigma_l.
		\end{equation*}
		 Using \eqref{sigma_bound}, we obtain
		\begin{align*}
			\sqrt{\delta}\left(\frac{1}{4C_*}-\frac{\sigma}{2}\right)&\le \frac{\sigma_l}{2} \frac{\sqrt{\delta}}{1-2 \sqrt{\delta} \sigma_l^2 \alpha_l}- \frac{\sigma }{2} \sqrt{\delta} 
			\\
			&\le \frac{\sqrt{\delta}}{2} \left( \sigma_l \left(\frac{1}{1-2 \sqrt{\delta} \sigma_l^2 \alpha_l}-1\right) + (\sigma_l-\sigma) \right) \\
			&\le \delta\sigma_l^3 \alpha_l + C  \delta^{3/2} .
		\end{align*}
	
		Therefore, we have
		\begin{equation}\label{B1G2}
			B_1-G_2 \le  \sigma_l^3 \alpha_l \delta \int_{y_0}^1 |f|^2 \, dy +  C\delta^{3/2} \int_{y_0}^1 |f|^2 \, dy. 
		\end{equation}
		
		\noindent$\bullet$ (Estimate of $D$):
		First, using $a \geq 1$ and change of variables, we estimate $D$ in terms of $f$ as
		\begin{align*}
			D &\geq \int_{\mathbb{R}_+} \frac{1}{v}|\partial_x (u-\tu)|^2 dx=\int_{y_0}^1 |\partial_y f|^2 \frac{1}{v} \left(\frac{dy}{dx} \right)\,dy.
		\end{align*}
		To estimate $\frac{dy}{dx}$ above, we use the following estimate: 
		 there exists $C>0$ such that
		\begin{equation}\label{Diffusion}
			\left|   \frac{1}{y(1-y)} \frac{1}{\tv} \left( \frac{dy}{dx} \right)-\frac{\sigma}{2 \sigma_l}\frac{\delta v''(p_-)}{ |v'(p_-)|^2 } \right| \le C \delta^2, \quad \forall x\in\bbr,
		\end{equation}
		where $v(p) = p^{-1/\gamma}$. We refer to \cite[Appendix A]{HKKL_pre}  for the proof of \eqref{Diffusion}. On the other hand, since $C^{-1} \leq v \leq C$, we have
		\begin{equation}\label{estimate in diffusion}
			\left| \frac{\tv}{v}-1 \right| \leq C \left|\tv-v \right| 
			\leq C \e.
		\end{equation}
		Then, using \eqref{Diffusion} and \eqref{estimate in diffusion}, we obtain the lower bound for $D$ as
		\begin{align*}
			D &\geq \int_{y_0}^1 |\partial_y f|^2 \frac{\tv}{v} \frac{1}{\tv} \left(\frac{dy}{dx} \right) dy \\
			&\geq (1-C\e)\left(\frac{\sigma}{2 \sigma_l}\frac{\delta v''(p_-)}{ |v'(p_-)|^2 }-C \delta^2 \right) \int_{y_0}^1 y(1-y) \left| \partial_y f \right|^2\,dy.
		\end{align*}
		Finally, since
		\begin{align*}
			\frac{1}{2 }\frac{v''(p_-)}{|v'(p_-)|^2 }=\frac{1}{2} (1+\gamma) \frac{1}{v_-}=\sigma_l^3 \alpha_l,
		\end{align*}
		we obtain
		\begin{align}
			\begin{aligned}\label{D_est} 
				D &\ge \sigma_l^3 \alpha_l \delta (1-C (\delta + \e)) \int_{y_0}^1 y(1-y) |\partial_y f|^2 \, dy \\   
				&\ge  \sigma_l^3 \alpha_l \delta (1-C (\delta + \e)) \int_{y_0}^1 (y-y_0)(1-y) |\partial_y f|^2 \, dy. 
			\end{aligned}    
		\end{align} 
		We combine the estimates \eqref{B1G2} and \eqref{D_est} to obtain
		\begin{align*}
		B_1-G_2-\frac{3}{4}D &\leq \sigma_l^3\alpha_l \delta\int_{y_0}^{1}|f|^2\,dy+C\delta^{3/2}\int_{y_0}^{1}|f|^2\,dy \\ 
			& \quad \ -\frac{3}{4}\sigma_l^3\alpha_l\delta(1-C(\delta+\e))\int_{y_0}^1(y-y_0)(1-y)|\pa_yf|^2\,dy.
		\end{align*}
		Now, we use Poincar\'e-type inequality in Lemma \ref{Poincare} with $c=y_0, d=1$ to have
		\[\bar{f}:=\frac{1}{1-y_0}\int_{y_0}^1f\,dy,\quad \mbox{and}\quad 
		\int_{y_0}^1 |f-\bar{f}|^2\,dy=\int_{y_0}^1 |f|^2\,dy-(1-y_0)\bar{f}^2.\]
		Thus,
		\begin{align}
		\begin{aligned}\label{est_B1G2D}
			B_1-G_2-\frac{3}{4}D &\leq \sigma_l^3\alpha_l \delta\int_{y_0}^{1}|f|^2\,dy+C\delta^{3/2}\int_{y_0}^{1}|f|^2\,dy \\ 
			& \quad \ -\frac{3}{2}\sigma_l^3\alpha_l\delta(1-C(\delta+\e))\left[\int_{y_0}^{1}|f|^2\,dy-(1-y_0)\bar{f}^2\right].
		\end{aligned}
		\end{align}
		Thus, we combine \eqref{xdot} and \eqref{est_B1G2D}, and using the smallness of $\delta$, $\e$, and $y_0$, we have
		\begin{align*}
			-\frac{\delta}{2M}&|\sX|^2+B_1-G_2-\frac{3}{4}D\\
			&\le -M\delta\left(\int_{y_0}^1f\,dy\right)^2+(\sigma_l^3\alpha_l\delta+C\delta^2+C\delta^{3/2})\int_{y_0}^1|f|^2\,dy\\
			&\quad -\frac{10}{8}\sigma_l^3\alpha_l\delta\left[\int_{y_0}^{1}|f|^2\,dy-(1-y_0)\bar{f}^2\right]\\
			&\le  -M\delta\left(\int_{y_0}^{1}f\,dy\right)^2+\frac{9}{8}\sigma_l^3\alpha_l\delta \int_{y_0}^{1}|f|^2\,dy   -\frac{10}{8}\sigma_l^3\alpha_l\delta\int_{y_0}^{1} |f|^2\,dy +\frac{5}{4}\frac{\sigma_l^3\alpha_l\delta}{1-y_0}\left(\int_{y_0}^{1}f\,dy\right)^2.
					\end{align*}
		Choosing $M=\frac{3}{2}\sigma_l^3\alpha_l$ and using the smallness $y_0<\frac{1}{6}$, we have 
		\[
		-\frac{\delta}{2M} |\sX|^2+B_1-G_2-\frac{3}{4}D \le -\frac{\sigma_l^3\alpha_l}{8}\delta\int_{y_0}^{1}|f|^2\,dy.
		\]
		 Therefore, taking $C_1=\frac{\sigma_l^3\alpha_l}{8}$, we get the desired inequality \eqref{est_R1}.
	\end{proof}

	\subsubsection{Estimate of $\mathcal{R}_2$}
	Now, we provide the estimate for the remaining terms $\mathcal{R}_2$. We use Cauchy-Schwarz inequality to estimate $\mathcal{R}_2$ as
	\begin{align*}
	\mathcal{R}_2&=-\frac{\delta}{2M}| {\dot{\bX}}|^2+\dot{\bX}\sum\limits_{i=3}^6\bY_i+\sum\limits_{i=2}^6B_i-G_1-\frac{D}{4}\\
	&\le-\frac{\delta}{4M}|\sX|^2 +\frac{C}{\delta}\sum_{i=3}^6 |\bY_i|^2 +\sum\limits_{i=2}^6B_i-G_1-\frac{D}{4}.
	\end{align*}
	We first obtain the following estimates on the bad terms in $\mathcal{R}_2$.
	
	\begin{lemma}\label{lem:est_R2}
		For sufficiently small $\delta$ and $\e$, 
		\[\frac{C}{\delta}\sum_{i=3}^6 |\bY_i|^2\le \frac{1}{8}(G_1+C_1G^S),\quad \sum_{i=2}^6B_i\le \frac{5}{32}(G_1+C_1G^S+D).\]
	\end{lemma}
	\begin{proof}
		Since the proof is essentially the same as in \cite[Section 4.5]{HKKL_pre}, we omit the proof.
	\end{proof}
	The estimates in Lemma \ref{lem:est_R2} yield the following bound on $\mathcal{R}_2$:
	\begin{align}
	\begin{aligned}\label{est_R2}
	\mathcal{R}_2&\le -\frac{\delta}{4M}|\sX|^2+\frac{9}{32}(G_1+C_1G^S)+\frac{5}{32}D-G_1-\frac{D}{4}\\
	&=-\frac{\delta}{4M}|\sX|^2+\frac{9}{32}C_1G^S-\frac{23}{32}G_1-\frac{3}{32}D.
	\end{aligned}
	\end{align}
	Therefore, from \eqref{decomml}, \eqref{est_R1}, and \eqref{est_R2}, we have
	\begin{align} \label{20}
		\frac{d}{dt}\int_{\mathbb{R}_+}a\eta(U | \widetilde{U})dx \le -\frac{\delta}{4M}|\sX|^2-\frac{23}{32}G_1-\frac{23}{32}C_1G^S-\frac{3}{32}D+\mathcal{P}.    
	\end{align}
	We now integrate \eqref{20} over $\left[0, t\right]$ for any $t\in \left[0, T\right]$ to derive
	
	\begin{align*}
	\int_{\R_+}&a(t,x)\eta(U(t,x)|\tU(t,x))\,dx+\int_0^t \left(\delta|\sX|^2+G_1+G^S+D\right)\,ds\\
	&\le C\left(\int_{\mathbb{R}_+} a(0,x)\eta(U(0,x)|\tU(0,x))\,dx+\int_0^t \mathcal{P}\,d s\right).
	\end{align*}
	Using $1\le a\le \frac{3}{2}$ and the estimate on the boundary term in Lemma \ref{lem:boundary}, we obtain
		\begin{align*}
			\int_{\mathbb{R}_+}&\eta(U(t,x)|\widetilde{U}(t,x))\,dx+\int_{0}^{t}\left(\delta|\sX|
			^2+G_1+G^S+D\right)\,ds \\
			\le& \ Ce^{-C\delta\beta}+C\int_{\mathbb{R}_+}\eta(U(0,x)|\widetilde{U}(0,x))dx+C\varepsilon^2\int_{0}^t\|(u-\widetilde{u})_{xx}\|_{L^2(\mathbb{R}_+)}^{2}ds.    
		\end{align*}
	Finally, by Lemma \ref{lem:relative}, the relative entropy is equivalent to the $L^2$-distance under the a priori assumption on small perturbation, that is,
	\[\|U(t,\cdot)-\tU(t,\cdot)\|_{L^2(\R_+)}^2\sim \int_{\R_+}\eta(U(t,x)|\tU(t,x))\,dx,\quad 0\le t\le T,\]
	we derive the desired estimate on the $L^2$-norm of the perturbation:
	\begin{align*}
		\lVert &(v-\widetilde{v})(t,\cdot) \rVert^2_{L^2(\mathbb{R}_+)}+\lVert (u-\widetilde{u})(t,\cdot) \rVert^2_{L^2(\mathbb{R}_+)}+\int_{0}^{t}\left(\delta|\sX|
		^2+G_1+G^S+D\right)\,ds\\
		&\le  C\left(\lVert (v-\widetilde{v})(0,\cdot) \rVert^2_{L^2(\mathbb{R}_+)}+\lVert (u-\widetilde{u})(0,\cdot) \rVert^2_{L^2(\mathbb{R}_+)}\right)+ Ce^{-C\delta\beta}+C\varepsilon^2\int_{0}^t\|(u-\widetilde{u})_{xx}\|_{L^2(\mathbb{R}_+)}^{2}ds.
	\end{align*}
	Finally, since $D\sim D_{u_1}$, we obtain the desired estimate \eqref{mainlemmaend}. This complete the proof of Lemma \ref{lem:main}.

	\section{$H^1$-perturbation estimate for the impermeable wall problem}\label{sec:5}
	\setcounter{equation}{0}
	In order to attain the $H^1$-estimate for the perturbation, we need another good term to control $\pa_x(v-\tv)$ term. Inspired by the previous literature \cite{KV21,KVW23} on the long-time behavior of the Navier-Stokes equations in the whole line, we introduce an effective velocity $h:=u-(\ln v)_x$ and consider the NS system in terms of $(v,h)$-variables:
	
	\begin{align}
	\begin{aligned} \label{heq1}
		&v_t-h_x=(\ln v)_{xx} \\
		&h_t+p(v)_x=0.
	\end{aligned}
	\end{align}
	The equations for the corresponding viscous shock wave in terms of $(\tv,\widetilde{h})$ become
	\begin{align*}
		&-\sigma\widetilde{v}'-\widetilde{h}'=(\ln \widetilde{v})'' \\
		&-\sigma\widetilde{h}'+(p(\widetilde{v}))'=0.
	\end{align*}
	Again, the system \eqref{heq1} can be written as general hyperbolic system of the form
	\[U_t +A(U)_x = \pa_x(M(U)\pa_x\nabla \eta(U)),\]
	where now the conserved quantity $U$, flux $A$ and the diffusion matrix $M$ are
	\[U=\begin{pmatrix}
	v\\ h
	\end{pmatrix},\quad A(U):=\begin{pmatrix}
	-h\\p(v)
	\end{pmatrix},\quad M(U):=\begin{pmatrix}
	\frac{1}{\gamma p(v)} & 0 \\ 0 & 0 
	\end{pmatrix}.\]
	Then, the entropy and relative entropy for this system are thus given by
	\begin{equation}\label{ent-vh}
	\eta(U) = \frac{|h|^2}{2}+Q(v),\quad \eta(U|\tU) = \frac{|h-\widetilde{h}|^2}{2}+Q(v|\tv).
	\end{equation}
	Similarly, the shifted shock wave  $(\tv,\widetilde{h}):= (\tv,\widetilde{h})(x-\sigma t -\bX(t) -\beta)$ for the system \eqref{heq1} satisfies
	\[
	\tU_t +A(\tU)_x = \pa_x(M(\tU)\pa_x\nabla\eta(\tU))-\sX\pa_x\tU.
	\]
	Then, using $(v,h)$ variables, we obtain the following estimate on the $H^1$-perturbation of $v$.
	\begin{lemma}\label{lem:v-H1}
		Under the hypothesis of Proposition \ref{prop:main}, there exists a positive constant $C$ such that for all $t\in[0,T]$,
		\begin{align}
		\begin{aligned}\label{vdiffusionlemmaend}
		\|v&-\tv\|_{H^1(\R_+)}^2+\|u-\tu\|_{L^2(\R_+)}^2 +\delta\int_0^t|\sX|^2\,ds+\int_0^t (G_1+G^S+D_{v_1}+D_{u_1})\,ds\\
		&\le C\left(\|v_0-\tv(0,\cdot)\|_{H^1(\R_+)}^2+\|u_0-\tu(0,\cdot)\|_{L^2(\R_+)}^2\right)+Ce^{-C\delta\beta}+C\varepsilon_1^2\int_{0}^t\|(u-\widetilde{u})_{xx}\|_{L^2(\mathbb{R}_+)}^{2}ds,
		\end{aligned}
		\end{align}
		where $D_{v_1}$ is the same term defined in \eqref{good_terms}.
	\end{lemma}
	\begin{proof}
	Following the proof of Lemma \ref{lem:rel_ent_est} with $a(t,x)\equiv 1$, we obtain the estimate on the time derivative of the relative entropy $\eta(U|\tU)$ given as \eqref{ent-vh}:
	\[{d \over dt}\int_{\mathbb{R}_+} \eta (U(t,x)|\tU(t,x)\, dx=\sX(t) \mathcal{Y}+\sum_{i=1}^5\mathcal{I}_i,\]
		where
		\begin{align}
		\begin{aligned}\label{est_v_h1}
		\mathcal{Y}&:=\int_{\mathbb{R}_+}\widetilde{h}_x(h-\widetilde{h})\,dx-\int_{\mathbb{R}_+}p'(\widetilde{v})\widetilde{v}_x(v-\widetilde{v})\,dx:=\mathcal{Y}_1+\mathcal{Y}_2,\\
		\mathcal{I}_1&:=-\int_{\mathbb{R}_+}\partial_x\left((p(v)-p(\widetilde{v}))(h-\widetilde{h})\right)\,dx,\quad \mathcal{I}_2:=-\int_{\mathbb{R}_+}\widetilde{h}_x p(v|\widetilde{v})\,dx, \\
		\mathcal{I}_3&:=\int_{\mathbb{R}_+}\left(p(v)-p(\widetilde{v})\right)\partial_x\left(\frac{1}{\gamma p(v)}\partial_x\left(p(v)-p(\widetilde{v})\right)\right)\,dx,\\
		\mathcal{I}_4&:=\int_{\mathbb{R}_+}\left(p(v)-p(\widetilde{v})\right)\partial_x\left(\frac{p(\widetilde{v})-p(v)}{\gamma p(v)p(\widetilde{v})}\partial_xp(\widetilde{v})\right)\,dx, \\
		\mathcal{I}_5&:=-\int_{\mathbb{R}_+}p(v|\widetilde{v})(\ln \widetilde{v})_{xx}\,dx.
		\end{aligned}
		\end{align}
		
		However, from the definition $h=u-(\ln v)_x$ and \eqref{shock_prop} we have
		\[|h-\widetilde{h}|\le C|u-\widetilde{u}|+C|\widetilde{v}_x||v-\widetilde{v}|+C|(v-\widetilde{v})_x|,\quad \mbox{and}\quad   |\widetilde{h}_x|\le C|\widetilde{v}_x|\le C|\tu_x|.\]
		Therefore, we control $\mathcal{Y}_1$ by using $|v-\tv|\le C|p(v)-p(\tv)|$ and Cauchy-Schwarz inequality as
		\begin{align*}
			|\mathcal{Y}_1|&\le C\int_{\R_+}|\tu_x||u-\tu|\,dx+C\int_{\R_+}|\tv_x||v-\tv|\,dx+C\int_{\R_+}|\tv_x||(v-\tv)_x|\,dx\\
			&\le C\int_{\mathbb{R}_+} (|\widetilde{u}_x|+|\tv_x|) |u-\widetilde{u}|\,dx +C\int_{\mathbb{R}_+} |\widetilde{v}_x| \left|p(v)-p(\widetilde{v})-\frac{u-\tu}{2C_*}\right|\,dx \\
			&\quad +C\int_{\mathbb{R}_+} |\widetilde{v}_x| |(v-\widetilde{v})_x|\,dx\\
			&\le C\sqrt{\delta}\sqrt{G^S}+C\delta\sqrt{G_1}+C\delta\lVert (v-\widetilde{v})_x \rVert_{L^2(\mathbb{R}+)}.
		\end{align*}
		Thus, we get
		\[C|\sX||\mathcal{Y}_1|\le \frac{\delta}{4}|\sX|^2+CG^S+C\delta G_1+C\delta\lVert (v-\widetilde{v})_x \rVert^2_{L^2(\mathbb{R}+)}.\]
		Similarly, we also obtain
		\[C|\sX||\mathcal{Y}_2|\le \frac{\delta}{4}|\sX|^2+CG^S+C\delta G_1.\]
		We now focus on estimating $\mathcal{I}_i$ for $i=1,2,\ldots,5$. We first consider the terms $\mathcal{I}_{1}$, $\mathcal{I}_3$, and $\mathcal{I}_4$. Using integration-by-parts, we get
		\begin{align*}
			\mathcal{I}_1&=\left[(p(v)-p(\widetilde{v}))(h-\widetilde{h})\right]_{x=0}, \\
			\mathcal{I}_3&=\left[-\frac{1}{\gamma p(v)}(p(v)-p(\widetilde{v}))\partial_x(p(v)-p(\widetilde{v}))\right]_{x=0}-\int_{\mathbb{R}_+}\frac{1}{\gamma p(v)} |(p(v)-p(\widetilde{v}))_x|^2\,dx \\
			&=:\mathcal{I}_{31}-\widetilde{D},\\
			\mathcal{I}_4&=-\left[(p(v)-p(\tv))^2\frac{\pa_x p(\tv)}{\gamma p(v)p(\tv)}\right]_{x=0}+\int_{\R_+}\pa_x(p(v)-p(\tv))\frac{p(v)-p(\tv)}{\gamma p(v)p(\tv)}\pa_xp(\tv)\,dx\\
			&=:\mathcal{I}_{41}+\mathcal{I}_{42}.
		\end{align*}		
		Substituting $h=u-(\ln v)_x$ and $p'(v)=-\gamma v^{-\gamma-1}$, we have
		\begin{align*}
			\mathcal{I}_1+\mathcal{I}_{31}=\left[(p(v)-p(\widetilde{v}))(u-\widetilde{u})\right]_{x=0}+\left[(p(v)-p(\widetilde{v}))\left(\frac{\widetilde{v}_x}{\widetilde{v}}-\frac{\widetilde{v}^{-\gamma-1}}{v^{-\gamma}}\widetilde{v}_x\right)\right]_{x=0}.   
		\end{align*}
		Since $u(t,0)=0$, we use the same argument as in the proof of Lemma \ref{lem:boundary} to obtain
		\begin{align*}
		\left[(p(v)-p(\widetilde{v}))(u-\widetilde{u})\right]_{x=0}\le \|p(v)-p(\tv)\|_{L^\infty(\R_+)}|\tu(t,0)|\le C\e\delta e^{-C\delta\beta}e^{-C\delta t},
		\end{align*} 
		and
		\begin{align*}
		\left[(p(v)-p(\widetilde{v}))\left(\frac{\widetilde{v}_x}{\widetilde{v}}-\frac{\widetilde{v}^{-\gamma-1}}{v^{-\gamma}}\widetilde{v}_x\right)\right]_{x=0}&\le \|p(v)-p(\tv)\|_{L^\infty(\R_+)}\left[\left|\frac{\tv_x}{\tv}\right|+\left|\frac{\tv^{-\gamma-1}}{v^{-\gamma}}\right||\tv_x|\right]_{x=0}\\
		&\le C\e \delta^2 e^{-C\delta\beta}e^{-C\delta t}.
		\end{align*}
		On the other hand, we estimate $\mathcal{I}_{41}$ as
		\[|\mathcal{I}_{41}|\le C\|p(v)-p(\tv)\|_{L^\infty(\R_+)}^2|\pa_xp(\tv(t,0))|\le C\e^2\delta^2e^{-C\delta\beta}e^{-C\delta t}.\]
		These estimates imply that 
		\[\mathcal{I}_1+\mathcal{I}_{31}+\mathcal{I}_{41}\le C\e\delta e^{-C\delta\beta}e^{-C\delta t}.\]
		Finally, we estimate $\mathcal{I}_{42}$ as
		\begin{align*}
		\mathcal{I}_{42}&\le \frac{1}{4}\int_{\R_+}\frac{|(p(v)-p(\tv))_x|^2}{\gamma p(v)}\,dx+C\int_{\R_+} |p(v)-p(\tv)|^2|\tv_x|^2\,dx \\
		&\le \frac{1}{4}\widetilde{D} +C\int_{\R_+}|\tu_x|^2\left|p(v)-p(\tv)-\frac{u-\tu}{2C_*}\right|^2\,dx+C\int_{\R_+}|\tu_x|^2|u-\tu|^2\,dx\\
		&\le  \frac{1}{4}\widetilde{D} +C\delta G_1+C\delta G^S.	
		\end{align*}
		Therefore, we obtain
		\begin{align*}
			\mathcal{I}_1+\mathcal{I}_3+\mathcal{I}_4&\le-\frac{3}{4}\widetilde{D}+C\delta G_1+C\delta G^S+ Ce^{-C\delta\beta}.\\
		\end{align*}
		To estimate $\mathcal{I}_2$, we first note that $\|v-\tv\|_{L^\infty(\R_+)}\le C\e$. Therefore, we can apply Lemma \ref{lem:relative} (3) to obtain
		\[0<p(v|\widetilde{v})\le C|p(v)-p(\widetilde{v})|^2.\]
		Since $|\widetilde{h}_x|\le C|\widetilde{v}_x|$ we have
		\[|\mathcal{I}_2|\le \sqrt{\delta}G_1+CG^S.\]
		Similarly, it is straightforward to estimate $\mathcal{I}_5$ as
		\begin{align*}
			|\mathcal{I}_5|&\le C\delta ^{3/2}G_1+CG^S.
		\end{align*}
		
		Thus, combining all the estimates above and integrating over $\left[0, t\right]$ for $0\le t\le T$, and using the fact that $\widetilde{D}\sim D_{v_1}$, we obtain
		\begin{align}
		\begin{aligned}\label{est-vh}
			\lVert &(v-\widetilde{v})(t, \cdot) \rVert^2_{L^2(\mathbb{R}_+)}+\lVert (h-\widetilde{h})(t, \cdot) \rVert^2_{L^2(\mathbb{R}_+)}+\int_{0}^{t}D_{v_1}\,ds\\
			&\le  C\left(\lVert (v-\widetilde{v})(0, \cdot) \rVert^2_{L^2(\mathbb{R}_+)}+\lVert (h-\widetilde{h})(0, \cdot) \rVert^2_{L^2(\mathbb{R}_+)}\right) +\frac{\delta}{2}\int_{0}^{t}|\sX(s)|^2\,ds\\
			&\quad+\sqrt{\delta}\int_{0}^{t}G_1\,ds+C_2\int_{0}^{t}G^S\,ds+C\e e^{-C\delta\beta},
		\end{aligned}
		\end{align}
		where $C_2>0$. Therefore, we multiply the estimates \eqref{est-vh} by $\frac{1}{2\max\left\{1,C_2 \right\}}$ and combine it with \eqref{mainlemmaend} to derive
		\begin{align*}
		&\|v-\tv\|^2_{L^2(\R_+)}+\|u-\tu\|^2_{L^2(\R_+)}+\|h-\widetilde{h}\|_{L^2(\R_+)}^2\\
		&\quad +\delta\int_0^t |\sX|^2\,ds+\int_0^t (G_1+G^S+D_{v_1}+D_{u_1})\,ds\\
		&\le C\left(\|v_0-\tv(0,\cdot)\|_{L^2(\R_+)}^2+\|u_0-\tu(0,\cdot)\|_{L^2(\R_+)}^2+\|h_0-\widetilde{h}(0,\cdot)\|_{L^2(\R_+)}^2\right)+Ce^{-C\delta\beta}\\
		&\quad +C\varepsilon^2\int_{0}^t\|(u-\widetilde{u})_{xx}\|_{L^2(\mathbb{R}_+)}^{2}ds.
		\end{align*}
		Finally, since 
		\[|(v-\widetilde{v})_x|\le C\left(|h-\widetilde{h}|+|v-\widetilde{v}|+|u-\widetilde{u}|\right)\quad\mbox{and}\quad |h-\widetilde{h}|\le C\left(|u-\widetilde{u}|+|v-\widetilde{v}|+|(v-\widetilde{v})_x|\right),\] 
		we obtain the desired $H^1$-estimate \eqref{vdiffusionlemmaend} for the perturbation on $v$.
	\end{proof}

	\subsection{Estimate of $\lVert (u-\widetilde{u})_x \rVert_{L^2}$} Finally, we complete the proof of the main proposition for the impermeable case by obtaining the $H^1$-norm of $u$ perturbation.
	
	\begin{lemma}
		Under the hypothesis of Proposition \ref{prop:main}, there exists a positive constant $C$ such that for all $t\in[0,T]$,

		\begin{align*}
		\|v&-\tv\|_{H^1(\R_+)}^2+\|u-\tu\|_{H^1(\R_+)}^2\\
		&\quad +\delta\int_0^t|\sX|^2\,ds+\int_0^t (G_1+G^S+D_{v_1}+D_{u_1}+D_{u_2})\,ds\\
		&\le C\left(\|v_0-\tv(0,\cdot)\|_{H^1(\R_+)}^2+\|u_0-\tu(0,\cdot)\|_{H^1(\R_+)}^2\right)+Ce^{-C\delta\beta},
		\end{align*}
		where $D_{u_2}$ is the term defined in \eqref{good_terms}.
	\end{lemma}
	
	\begin{proof}
		Recall that perturbations $v-\tv$ and $u-\tu$ satisfy
		\begin{align}
		\begin{aligned} \label{oeh1}
			&(v-\widetilde{v})_t=(u-\widetilde{u})_x+\sX(t)\widetilde{v}_x, \\
			&(u-\widetilde{u})_t=-(p(v)-p(\widetilde{v}))_x+\left(\frac{u_x}{v}-\frac{\widetilde{u}_x}{\widetilde{v}}\right)_x+\sX(t)\widetilde{u}_x.
		\end{aligned} 
		\end{align}
		We multiply \eqref{oeh1}$_2$ by $-(u-\widetilde{u})_{xx}$ and integrating over $\left[0, t\right]\times \mathbb{R}_+$ to get
		\begin{align}
			\begin{aligned} \label{F1F2F3}
				-\int_{0}^{t}&\int_{\R_+}(u-\widetilde{u})_t(u-\widetilde{u})_{xx}\,dxds\\
				&=\int_{0}^{t}\int_{\R_+} (p(v)-p(\widetilde{v}))_x(u-\widetilde{u})_{xx}\,dxds -\int_{0}^{t}\int_{\R_+}\left(\frac{u_x}{v}-\frac{\widetilde{u}_x}{\widetilde{v}}\right)_x(u-\widetilde{u})_{xx}\,dxds \\
				&\quad  -\int_{0}^{t}\int_{\R_+}\sX(s)\widetilde{u}_x(u-\widetilde{u})_{xx}\,dxds .
			\end{aligned}
		\end{align}
		Applying the integration-by-parts to the l.h.s., we have
		\begin{align*}
			-\int_{0}^{t}&\int_{\R_+}(u-\widetilde{u})_t(u-\widetilde{u})_{xx}\,dxds\\
			&=\int_{0}^{t}\left[(u-\widetilde{u})_t(u-\widetilde{u})_x \right]_{x=0}\,ds+\int_{0}^{t} \frac{d}{ds}\int_{\R_+}\frac{|(u-\widetilde{u})_x|^2}{2}\,dxds \\
			&=\int_{0}^{t}\left[(u-\widetilde{u})_t(u-\widetilde{u})_x \right]_{x=0}\,ds+\frac{1}{2}\lVert (u-\widetilde{u})_x(t,\cdot) \rVert^2_{L^2(\mathbb{R}_+)}-\frac{1}{2}\lVert (u-\widetilde{u})_x(0,\cdot) \rVert^2_{L^2(\mathbb{R}_+)}.
		\end{align*}
		Thus, we rewrite \eqref{F1F2F3} as 
		\begin{equation} \label{24}
			\frac{1}{2}\lVert (u-\widetilde{u})_x(t,\cdot) \rVert^2_{L^2(\mathbb{R}_+)}=\frac{1}{2}\lVert (u-\widetilde{u})_x(0,\cdot) \rVert^2_{L^2(\mathbb{R}_+)}+F_1+F_2+F_3+F_4,
		\end{equation}
		where 
		\begin{align*}
		&F_1:=\int_{0}^{t}\int_{\R_+} (p(v)-p(\widetilde{v}))_x(u-\widetilde{u})_{xx}\,dxds,\\
		&F_2:= -\int_{0}^{t}\int_{\R_+}\left(\frac{u_x}{v}-\frac{\widetilde{u}_x}{\widetilde{v}}\right)_x(u-\widetilde{u})_{xx}\,dxds,\\
		&F_3:= -\int_{0}^{t}\int_{\R_+}\sX(s)\widetilde{u}_x(u-\widetilde{u})_{xx}\,dxds,\\
		&F_4:=-\int_{0}^{t}\left[(u-\widetilde{u})_t(u-\widetilde{u})_x\right]_{x=0}\,ds.
		\end{align*}
		Below, we estimate $F_i$ for $i=1,\ldots,4$.\\

		\noindent $\bullet$ (Estimate of $F_1$): Using Young's inequality, we have
		\begin{align*}
			|F_1|&\le C\int_{0}^{t}\int_{\R_+} |(p(v)-p(\widetilde{v}))_x|^2 \,dx\,ds+\frac{1}{4}\int_0^t\int_{\R_+}\frac{|(u-\widetilde{u})_{xx}|^2}{v} \,dx\,ds.\\ 
		\end{align*}
		
		\noindent$\bullet$ (Estimate of $F_2$): To estimate $F_2$, we split it as
		\begin{align*}
			F_2&=-\int_{0}^{t}\int_{\R_+} \frac{|(u-\widetilde{u})_{xx}|^2}{v}\,dxds-\int_{0}^{t}\int_{\R_+} \left(\frac{1}{v}\right)_x(u-\widetilde{u})_x(u-\widetilde{u})_{xx}\,dxds \\
			&\quad-\int_{0}^{t}\int_{\R_+} \left(\widetilde{u}_x\left(\frac{1}{v}-\frac{1}{\widetilde{v}}\right)\right)_x\,dxds \\
			&=:-\int_{0}^{t}\int_{\R_+} \frac{|(u-\widetilde{u})_{xx}|^2}{v}\,dxds+F_{21}+F_{22}.
		\end{align*}
		Then, we again separate $F_{21}$ as
		\begin{align*}
			F_{21}&=\int_{0}^{t}\int_{\R_+} \frac{v_x}{v^2}(u-\widetilde{u})_x(u-\widetilde{u})_{xx}\,dxds\\    
			&=\int_{0}^{t}\int_{\R_+} \frac{(v-\widetilde{v})_x}{v^2}(u-\widetilde{u})_x(u-\widetilde{u})_{xx}\,dxds+\int_{0}^{t}\int_{\R_+} \frac{\widetilde{v}_x}{v^2}(u-\widetilde{u})_x(u-\widetilde{u})_{xx}\,dxds\\
			&=:F_{211}+F_{212}.
		\end{align*}
		We use the a priori assumption, Holder's inequality and the Sobolev interpolation inequality to obtain
		\begin{align*}
			|F_{211}|&\le C\int_{0}^{t} \lVert (v-\widetilde{v})_x \rVert_{L^2(\mathbb{R}+)}\lVert (u-\widetilde{u})_x \rVert_{L^\infty(\mathbb{R}+)}\lVert (u-\widetilde{u})_{xx} \rVert_{L^2(\mathbb{R}+)} \,ds \\
			&\le C\e\int_{0}^{t} \lVert (u-\widetilde{u})_x \rVert^{\frac{1}{2}}_{L^2(\mathbb{R}+)} \lVert (u-\widetilde{u})_{xx} \rVert^{\frac{3}{2}}_{L^2(\mathbb{R}+)}\,ds   \\
			& \le C\e\int_{0}^{t} \left(\lVert (u-\widetilde{u})_{x} \rVert^2_{L^2(\mathbb{R}+)}+\lVert (u-\widetilde{u})_{xx} \rVert^2_{L^2(\mathbb{R}+)}\right)\,ds \\
			&\le C\e\int_0^t\left(D_{u_1}+D_{u_2}\right)\,ds.
		\end{align*}
		We estimate $F_{212}$ in a similar manner as
		\begin{align*}
			|F_{212}|\le C\delta^2\int_0^t\left(D_{u_1}+D_{u_2}\right)\,ds.    
		\end{align*}
		To estimate $F_{22}$, we use Young's inequality to derive
		\begin{align*}
			\left|F_{22}\right|&=\left|\int_{0}^{t}\int_{\R_+} \left[\frac{\widetilde{u}_{xx}}{v\widetilde{v}}(v-\widetilde{v})+\frac{\widetilde{u}_x}{v^2}(v-\widetilde{v})_x+\frac{\widetilde{u}_x\widetilde{v}_x}{v^2\widetilde{v}^2}(v^2-\widetilde{v}^2)\right]\,dxds\right| \\
			&\le C\delta^2\int_{0}^{t}\int_{\R_+} \left(|v-\widetilde{v}|^2+|(v-\widetilde{v})_x|^2+|(u-\widetilde{u})_{xx}|^2\right)\,dxds \\
			&\le C\delta^2\int_{0}^{t} \left(G_1+G^S+D_{v_1}+D_{u_2}\right)\,ds.
		\end{align*}
		\noindent $\bullet$ (Estimate of $F_3$): We use Cauchy-Schwarz inequality to obtain
		\begin{align*}
			\left|F_3\right|&\le C\delta^2 \int_{0}^{t} \int_{\R_+} |\sX(s)| |(u-\widetilde{u})_{xx}|\,dx ds \\ 
   &\le C\delta^2\int_{0}^{t} |\sX(s)|^2\, ds+\frac{1}{8}\int_{0}^{t}\int_{\mathbb{R}_+} \frac{|(u-\widetilde{u})_{xx}|^2}{v}\,dxds. 
		\end{align*}
		
		\noindent $\bullet$ (Estimate of $F_4$): Finally, to estimate the boundary term $F_4$, we use \eqref{oeh1}$_1$ to observe that
		\begin{align*}
		F_4&=-\int_{0}^{t}\left[(u-\widetilde{u})_t(v-\widetilde{v})_t\right]_{x=0}\,ds+\int_{0}^{t}\sX(s)\left[(u-\widetilde{u})_t\widetilde{v}_x\right]_{x=0}\,ds\\
		&=(u-\tu)_t(0,0)(v-\tv)(0,0)-(u-\tu)_t(t,0)(v-\tv)(t,0)+\int_{0}^{t}\left[(u-\widetilde{u})_{tt}(v-\widetilde{v})\right]_{x=0}\,ds\\
		&\quad +\int_{0}^{t}\sX(s)\left[(u-\widetilde{u})_t\widetilde{v}_x\right]_{x=0}\,ds.
		\end{align*}
		Note that the impermeable boundary condition $u(t,0)=0$ implies $u_t(t,0)=u_{tt}(t,0)=0$. Moreover, 
	using the same argument as before, together with 
	\[
	|(\widetilde{u}_t)_{x=0}|, |(\tu_{tt})_{x=0}|\le C |\tu'(-\sigma t -\bX(t)-\beta)| |\sigma + \dot{\bX}(t)| \le C\delta^2 e^{-C(\sigma t+\beta)\delta},\quad t\le T,
	\]
	we have
		\[|F_4|\le C\delta\e e^{-C\delta\beta}.\]
		
		Therefore, combining the estimates for $F_i$ for $i=1,2,3,4$ and substituting them to \eqref{24}, we conclude that
		\begin{align*}
			\lVert &(u-\widetilde{u})_x(t,\cdot) \rVert^2_{L^2(\mathbb{R}_+)}+\frac{5}{8}\int_{0}^{t}\int_{\R_+}\frac{|(u-\tu)_{xx}|^2}{v}\,dx\,ds\\
			&\le \lVert (u-\widetilde{u})_x(0,\cdot) \rVert^2_{L^2(\mathbb{R}_+)}+C\delta\e e^{-C\delta\beta}+C\int_{0}^{t} \left(\delta|\sX|^2+G_1+G^S+D_{v_1}+D_{u_1}\right)\,ds \\
			&\quad +C(\e+\delta^2)\int_0^tD_{u_2}\,ds.
		\end{align*}
	However, since
	\begin{equation}\label{Du2_equiv}
	\int_{\R_+}\frac{|(u-\tu)_{xx}|^2}{v}\,dx\sim D_{u_2},
	\end{equation}
	we have
	\begin{align}
	\begin{aligned} \label{25}
	\lVert &(u-\widetilde{u})_x(t,\cdot) \rVert^2_{L^2(\mathbb{R}_+)}+\frac{1}{2}\int_{0}^{t}\int_{\R_+}\frac{|(u-\tu)_{xx}|^2}{v}\,dx\,ds\\
	&\le \lVert (u-\widetilde{u})_x(0,\cdot) \rVert^2_{L^2(\mathbb{R}_+)}+C\delta\e e^{-C\delta\beta}+C_3\int_{0}^{t} \left(\delta|\sX|^2+G_1+G^S+D_{v_1}+D_{u_1}\right)\,ds.
	\end{aligned}
	\end{align}

	We now multiply \eqref{25} by $\frac{1}{2\max\left\{1, C_3  \right\}}$ and then add to \eqref{vdiffusionlemmaend} and then again use \eqref{Du2_equiv} to derive
	\begin{align*}
		\lVert (v-&\widetilde{v})(t,\cdot) \rVert^2_{H^1(\mathbb{R}_+)}+\lVert (u-\widetilde{u})(t,\cdot) \rVert^2_{H^1(\mathbb{R}_+)} \\ 
		&+\int_{0}^{t}\left(\delta|\sX|
		^2+G_1+G^S+D_{v_1}+D_{u_1}+D_{u_2}\right)ds \\  
		\le &  C\left(\lVert (v-\widetilde{v})(0, \cdot) \rVert^2_{H^1(\mathbb{R}_+)}+\lVert (u-\widetilde{u})(0, \cdot) \rVert^2_{H^1(\mathbb{R}_+)}\right)+Ce^{-C\delta\beta},
	\end{align*}
	which is the desired estimate on the $H^1$-perturbation. This complete the proof of Proposition \ref{prop:main}.
	
	\end{proof}
	
	\section{Perturbation estimates for the inflow problem}\label{sec:6}
	\setcounter{equation}{0}
	
	In this section, we consider the inflow problem \eqref{eq:NS_inflow_2}--\eqref{boundary_inflow}. Below, we provide $L^2$ and $H^1$-perturbation estimates and thereby, prove Proposition \ref{prop:main} for the inflow problem. Since the proofs share a lot of parts with the proof of Proposition \ref{prop:main} in the previous sections, we focus on delivering the difference between the inflow problem and the impermeable wall problem.
	
	\subsection{Relative entropy estimate for the inflow problem}
	
	Let us first present the estimate on the $L^2$-perturbation for the inflow problem.
	
	\begin{lemma}\label{lem:main-2}
		Under the hypothesis of Proposition \ref{prop:main}, there exists a positive constant $C$ such that for all $t\in[0,T]$,
		\begin{align}
		\begin{aligned} \label{mainlemmaend-2}
		&\|v-\tv\|_{L^2(\R_+)}^2+\|u-\tu\|_{L^2(\R_+)}^2+\int_0^t (\delta|\sX(s)|^2+G_1+G^S+D_{u_1})\,ds\\
		&\quad \le C\|v_0-\tv(0,\cdot)\|_{L^2(\R_+)}^2+\|u_0-\tu(0,\cdot)\|^2_{L^2(\R_+)}+Ce^{-C\delta \beta}+ C\varepsilon^{2}\int_{0}^{t}\|(u-\widetilde{u})_{\xi\xi}\|^{2}_{L^2}\,ds,
		\end{aligned}
		\end{align}
		where $G_1$, $G^S$, and $D_{u_1}$ are the terms defined in \eqref{good_terms}.
	\end{lemma}

%
	Since the inflow problem is written in terms of $\xi=x-\sigma_-t$, we consider the weight function $a$ in $\xi$-variable:
	\begin{equation}\label{weight_inflow}
	a(t,\xi):=1+\frac{u_--\tu(t,\xi)}{\sqrt{\delta}}=1+\frac{u_--\tu(\xi-(\sigma-\sigma_-)t-\bX(t)-\beta)}{\sqrt{\delta}}.
	\end{equation}
	Then, the weight function $a$ still satisfies $1\le a<\frac{3}{2}$ and 
	\[\pa_\xi a = -\frac{\pa_\xi \tu}{\sqrt{\delta}}=\frac{\sigma \pa_\xi \tv}{\sqrt{\delta}}>0,\quad \mbox{and}\quad |a_\xi|\sim \frac{\pa_\xi\tv}{\sqrt{\delta}}.\]
	Moreover, the system \eqref{eq:NS_inflow_2} can be written in the same structure as before
	\begin{equation} \label{general_inflow}
	\rd_t U+ \rd_{\xi} A(U)=\rd_{\xi} \left( M(U) \rd_{\xi} \nabla \eta (U) \right),
	\end{equation}
	where
	\[U=\begin{pmatrix}
	v\\u
	\end{pmatrix},\quad M(U) = \begin{pmatrix}
	0 & 0 \\ 0 & \frac{1}{v}
	\end{pmatrix},\quad \eta(U) = \frac{u^2}{2}+Q(v),\] 
	and the flux $A$ now becomes
	\[
	A(U):=\begin{pmatrix}
	-\sigma_-v-u\\
	-\sigma_-u+p(v)
	\end{pmatrix}.\]
	Then, the relative quantities for the system \eqref{general_inflow} can be written as
	\begin{equation} \label{relative_quantities_inflow} 
	\begin{split}
	&\eta(U|\tU)={|u-\tu|^2 \over 2}+Q(v|\tv),\quad A(U|\tU)=\begin{pmatrix}
	0\\
	p(v|\tv)
	\end{pmatrix},\\
	&G(U;\tU)=(p(v)-p(\tv))(u-\tu) - \sigma_-\eta(U|\tU).
	\end{split}
	\end{equation}
	
	Moreover, the shifted viscous shock profile $\tU$ 
	\[
	\tU(t,\xi)=\begin{pmatrix}
	\tv(t,\xi)\\
	\tu(t,\xi)
	\end{pmatrix}=
	\begin{pmatrix}
	\tv(\xi-(\sigma-\sigma_-) t-\bX(t)-\beta)\\
	\widetilde{u}(\xi-(\sigma-\sigma_-) t-\bX(t)-\beta)
	\end{pmatrix}
	\]
	satisfies
	\begin{equation} \label{viscous_inflow} 
	\rd_t \tU+\rd_{\xi} A(\tU)=\rd_{\xi} \left( M(\tU) \rd_{\xi} \nabla \eta (\tU) \right)-\sX \rd_{\xi}\tU.
	\end{equation}
	
	Below, we estimate the weighted relative entropy between the solution $U$ of \eqref{general_inflow} and the shifted viscous shock wave \eqref{viscous_inflow}. 
	\begin{lemma} 
		Let $a$ be the weight function defined by \eqref{weight_inflow}, $U$ be a solution to \eqref{general_inflow},  and $\tU$ be the shifted viscous shock wave satisfying \eqref{viscous_inflow}. Then,
		\[
		{d \over dt}\int_{\mathbb{R}_+} a(t,\xi)\eta (U(t,x))|\tU(t,\xi)) d \xi
		=\sX(t)\bold{Y}(U)+\mathcal{J}^{\textup{bad}}(U)-\mathcal{J}^{\textup{good}}(U)+\mathcal{P}(U),
		\]
		where
		\begin{align}
		\begin{aligned}\label{terms_rel_ent_inflow}
		\bold{Y}(U)&:= -\int_{\mathbb{R}_+} a_{\xi} \eta (U|\tU)\,d\xi  +\int_{\mathbb{R}_+} a\nabla^2 \eta (\tU)\tU_{\xi}(U-\tU)\,d\xi   , \\
		\mathcal{J}^{\textup{bad}}(U)&:=\int_{\mathbb{R}_+} a_{\xi}(p(v)-p(\tv))(u-\tu)\,d\xi  -\int_{\mathbb{R}_+} a \widetilde{u}_{\xi}p(v | \tv) \,d \xi \\
		&\quad-\int_{\mathbb{R}_+} a_{\xi} {u-\tu \over v} \rd_{\xi} (u-\tu) \,d\xi  + \int_\R a_{\xi}(u-\tu)(v-\tv){\tu_{\xi} \over v \tv} \,d\xi  \\
		&\quad +\int_{\mathbb{R}_+} a\rd_{\xi} (u-\tu) {v-\tv \over v \tv} \tu_{\xi} \,d\xi,   \\
		\mathcal{J}^{\textup{good}}(U)&:={\sigma \over 2}\int_{\mathbb{R}_+} a_{\xi}|u-\tu|^2\,d\xi  +\sigma \int_\R a_{\xi} Q(v|\tv)\,d\xi + \int_{\mathbb{R}_+} {a \over v} |\rd_{\xi} (u-\tu)|^2\,d\xi  ,\\
		\mathcal{P}(U)&:=\Big[a(u-\widetilde{u})(p(v)-p(\widetilde{v})) -a\frac{\sigma_-}{2}(u-\tu)^{2} - a\sigma_-Q(v|\tv) \\
		& \hspace{2.5cm} -\frac{a}{v}(u-\widetilde{u})(u-\widetilde{u})_{\xi}-\frac{a}{v\tv}(u-\tu)(v-\tv)\tu_{\xi}\Big]_{\xi=0}.
		\end{aligned}
		\end{align}
	\end{lemma}
	
	\begin{proof}
		Since the proof is almost the same as the proof of Lemma \ref{lem:rel_ent_est}, we only explain the different points. Again, the time derivative of the weighted relative entropy is given as
		\begin{equation}\label{rel_ent_est_inflow}
		{d \over dt}\int_{\mathbb{R}_+} a(t,\xi) \eta \left(U(t,\xi)|\tU(t,\xi)\right) d\xi=\sX(t) \bY -(\sigma-\sigma_-)\int_{\mathbb{R}_+}a_\xi\eta(U|\widetilde{U})d\xi+\sum_{i=1}^5I_i,
		\end{equation}
		where $I_i$ are given as in \eqref{I_i} in terms of $\xi$-variable. Note that, the coefficient in front of the relative entropy on the right-hand side of \eqref{rel_ent_est_inflow} becomes $(\sigma-\sigma_-)$, instead of $\sigma$. Then, compared to the estimates in the proof of Lemma \ref{lem:rel_ent_est}, the only difference is that $G(U;\widetilde{U})$ is changed from \eqref{relative_quantities} to \eqref{relative_quantities_inflow}. Due to this change, the term $I_1$ becomes
		\begin{align*}
		I_1&=a\left[(u-\tu)(p(v)-p(\tv))-\frac{\sigma_-}{2}(u-\tu)^2-\sigma_-Q(v|\tv)\right]_{\xi=0}\\
		&\quad +\int_{\R_+}a_\xi(u-\tu)(p(v)-p(\tv))\,d\xi-\frac{\sigma_-}{2}\int_{\R_+}a_\xi (u-\tu)^2\,d\xi-\sigma_-\int_{\R_+}a_\xi Q(v|\tv)\,d\xi,
		\end{align*}
		and the other terms $I_i$ for $i=2,3,4,5$ are unchanged, except that the boundary value of $u(t,0)$ does not vanish in the case of inflow problem. This induces that the boundary terms $\mathcal{P}(U)$ are slightly changed compared to that of the impermeable wall problem. Thus, combining all the terms $I_i$ for $i=1,2,\ldots, 5$, we get the desired estimate on the weighted relative entropy.
	\end{proof} 
	
	Observe that, all the terms in \eqref{terms_rel_ent_inflow} are the same as in \eqref{terms_YJP} and the only difference appears in the boundary terms $\mathcal{P}(U)$. In the following lemma, we provide the estimates on the boundary terms for the case of the inflow problem.
	
	\begin{lemma} \label{lem:boundary_inflow}
		There exists $C>0$ such that
		\[
		\left|\int_{0}^t\mathcal{P}(U) \,ds\right|  \le Ce^{-C\delta \beta}  + C\varepsilon^{2}\int_{0}^{t}\|(u-\widetilde{u})_{\xi\xi}\|_{L^2(\mathbb{R}+)}^{2}\,ds
		\]
		for all $t\in \left[0, T\right]$.
	\end{lemma}
	\begin{proof}
		We split $\mathcal{P}(U)$ as
		\begin{align*}
		\mathcal{P}(U)&:=\big[a(u-\widetilde{u})(p(v)-p(\widetilde{v})) -a\frac{\sigma_-}{2}(u-\tu)^{2} - a\sigma_-Q(v|\tv) \\
		& \quad -\frac{a}{v}(u-\widetilde{u})(u-\widetilde{u})_{\xi}-\frac{a}{v\tv}(u-\tu)(v-\tv)\tu_{\xi}\big]_{\xi=0} \\
		&=:\mathcal{P}_1(U)+\mathcal{P}_2(U)+\mathcal{P}_3(U)+\mathcal{P}_4(U)+\mathcal{P}_5(U).
		\end{align*}     
		First of all, since $\sigma>0$, $\sigma_-<0$, and $\beta>0$, we have \[|-(\sigma-\sigma_-) t-\bX(t)-\beta|\geq(\sigma-\sigma_-) t-|\bX(t)|+\beta,\]
		which together with $|\bX(t)|\leq C\varepsilon t$, we have
		\[|-(\sigma-\sigma_-) t-\bX(t)-\beta|\ge\frac{1}{2}(\sigma-\sigma_-)t+\beta.\] 
		This yields
		\begin{equation} \label{betaest-2}
		\begin{aligned}
		&|\widetilde{u}(t,0)-u_{-}|\le C\delta e^{-C\delta|-(\sigma-\sigma_-) t-\bX(t)-\beta|}=C\delta e^{-C\delta (\sigma-\sigma_-)t}e^{-C\delta \beta}, \\
		&|\widetilde{v}(t,0)-v_{-}|\le C\delta e^{-C\delta|-(\sigma-\sigma_-) t-\bX(t)-\beta|}=C\delta e^{-C\delta (\sigma-\sigma_-)t}e^{-C\delta \beta}.
		\end{aligned}
		\end{equation}
		Below, we estimate $\mathcal{P}_i(U)$ separately. First, we use a priori assumption and \eqref{betaest-2} to obtain
		\begin{align*}
		\left|\int_{0}^{t}\mathcal{P}_1(U)\,ds\right|\le C\lVert p(v)-p(\widetilde{v}) \rVert_{L^\infty(\mathbb{R}+)} \int_{0}^{t} |\widetilde{u}(s,0)-u_-| \,ds \le C\varepsilon e^{-C\delta\beta}.    
		\end{align*}
		Similarly, we estimate $\mathcal{P}_2(U)$ and $\mathcal{P}_3(U)$ as 
		\[\left|\int_{0}^{t}\mathcal{P}_2(U)\,ds\right|\le C\lVert u-\widetilde{u} \rVert_{L^\infty(\mathbb{R}+)} \int_{0}^{t} |\widetilde{u}(s,0)-u_-|\,ds \le C\varepsilon e^{-C\delta\beta},\]
		and
		\[\left|\int_{0}^{t}\mathcal{P}_3(U)\,ds\right|\le C\lVert v-\widetilde{v} \rVert_{L^\infty(\mathbb{R}+)} \int_{0}^{t} |\widetilde{v}(s,0)-v_-|\,ds \le C\varepsilon e^{-C\delta\beta}.\]
		
		For $\mathcal{P}_4(U)$, we use interpolation inequality and Young's inequality to get
		
		\begin{align*}
		\left|\int_{0}^{t}\mathcal{P}_4(U)\,ds\right| &= \left|\int_{0}^{t}\left[\frac{a}{v}(u-\widetilde{u})(u-\widetilde{u})_{\xi}\right]_{\xi=0} \,ds \right| \\
		& \leq C \int_{0}^{t} |u_- -\tu(s,0)| \|(u-\widetilde{u})_{\xi}\|_{L^\infty(\mathbb{R}+)}\,ds \\
		& \leq C\int_{0}^{t}  |u_- -\tu(s,0)|^{\frac{4}{3}}\,ds + C\int_{0}^{t} \|(u-\widetilde{u})_{\xi}\|^{2}_{L^2(\mathbb{R}+)}\|(u-\widetilde{u})_{\xi\xi}\|^{2}_{L^2(\mathbb{R}+)}\,ds \\
		& \leq C\delta^{\frac{1}{3}}e^{-C\delta\beta} + C\varepsilon^{2}\int_{0}^{t}\|(u-\widetilde{u})_{\xi\xi}\|^{2}_{L^2(\mathbb{R}+)}\,ds.
		\end{align*}
		 Similarly, we also easily get the estimate of $\mathcal{P}_5(U)$ as
		\[\left|\int_{0}^{t}\mathcal{P}_{5}(U)\,ds\right|\le C\varepsilon \delta^2e^{-C\delta\beta}.\]
		
		Combining the estimates on $\mathcal{P}_i(U)$ above with smallness of parameters, we obtain the desired estimate.
	\end{proof}
	
	With the above control on the boundary terms, the remaining process for proving Lemma \ref{lem:main-2} is the same as the proof of Lemma \ref{lem:main} in Section \ref{sec:proof_main_1}, except the small changes that we need to rewrite the $x$-variables as $\xi+\sigma_-t$. For example, when we obtain the counterpart of Lemma \ref{lem:est_R1}, we need to define the $f$ and $y$ variables as
	\[
	f:=u-\tu(\xi-(\sigma-\sigma_-) t-\bX(t)-\beta),
	\]
	and
	\[
	y:=\frac{u_--\tu(\xi-(\sigma-\sigma_-) t-\bX(t)-\beta)}{\delta},
	\]
	so that these are written in terms of $\xi$ and $t$. Except for replacing $x$ by $\xi$, the entire procedure is the same, and we conclude that the desired estimate \eqref{mainlemmaend-2} in Lemma \ref{lem:main-2} holds.
	
	\subsection{$H^1$-estimate for the inflow problem}
	To attain the $H^1$-estimate for the inflow problem, we again introduce an effective velocity $h:=u-(\ln v)_{\xi}$. Then, the NS system in terms of $(v,h)$ variable becomes
	\begin{align*}
	&v_t-\sigma_- v_{\xi} - h_{\xi}=(\ln v)_{\xi\xi}, \\
	&h_t- \sigma_-h_{\xi}+p(v)_{\xi}=0,
	\end{align*}
	and the equations for corresponding viscous shock wave becomes
	\begin{align*}
	&\widetilde{v}_t-\sigma_-\widetilde{v}_{\xi}-\widetilde{h}_{\xi}=(\ln \widetilde{v})_{\xi\xi} - \sX\widetilde{v}_{\xi}, \\
	&\widetilde{h}_t-\sigma_-\widetilde{h}_{\xi}+p(\widetilde{v})_{\xi}=-\sX\widetilde{h}_{\xi},
	\end{align*}
	where $\widetilde{h} := \widetilde{u}-(\ln\widetilde{v})_{\xi}$. Similar to Lemma \ref{lem:v-H1}, we can obtain the following $H^1$-estimate for $v$-perturbation.
	
	\begin{lemma}\label{fistorderv}
		Under the hypothesis of Proposition \ref{prop:main}, there exists a positive constant $C$ such that for all $t\in[0,T]$,
		\begin{align}
		\begin{aligned} \label{firstorderv}
		&\lVert (v-\widetilde{v})(t, \cdot) \rVert^2_{H^1(\mathbb{R}_+)}+\lVert (u-\widetilde{u}(t, \cdot) \rVert^2_{L^2(\mathbb{R}_+)} +\delta \int_{0}^{t}|\sX|^2\,ds+\int_{0}^{t}\left(G_1+G^S+D_{v_1} +D_{u_1}\right)\,ds\\   
		& \le  C\left(\lVert (v-\widetilde{v})(0, \cdot) \rVert^2_{H^1(\mathbb{R}_+)}+\lVert (u-\widetilde{u})(0, \cdot) \rVert^2_{L^2(\mathbb{R}_+)}\right)
		+ Ce^{-C\delta\beta} + \frac{1}{100}\int_{0}^{t}D_{u_2}\,ds,
		\end{aligned}
		\end{align}
		where $D_{v_1}, D_{u_2}$ are the term defined in \eqref{good_terms}.
	\end{lemma}
	\begin{proof}
		
		As in the proof of Lemma \ref{lem:v-H1}, we can estimate the relative entropy between $U=(v,h)$ and $\tU:=(\tv,\widetilde{h})$ as
		
		\[{d \over dt}\int_{\mathbb{R}_+} \eta (U(t,\xi)|\tU(t,\xi)) d\xi  =\sX(t) \mathcal{Y}(U)+\sum_{i=1}^5\mathcal{I}_i,\]
		where
		\begin{align*}
		&\mathcal{Y}(U):= \int_{\mathbb{R}_+} \nabla^2 \eta (\tU)(\tU)_{\xi}(U-\tU)\,d\xi, \\
		&\mathcal{I}_1:=-\int_{\mathbb{R}_+}  \rd_{\xi} G(U; \tU) \,d\xi,\\
		&\mathcal{I}_2:=-\int_{\mathbb{R}_+}  \rd_{\xi} \nabla \eta (\tU) A(U|\tU)\,d \xi,\\
		&\mathcal{I}_3:=\int_{\mathbb{R}_+}  \left( \nabla \eta (U)-\nabla \eta(\tU)\right) \rd_{\xi} \left(M(U) \rd_{\xi}(\nabla \eta(U)-\nabla \eta(\tU) \right)\,d\xi,\\
		&\mathcal{I}_4:= \int_{\mathbb{R}_+} \left( \nabla \eta(U)-\nabla \eta (\tU)\right) \rd_{\xi} \left( (M(U)-M(\tU))\rd_{\xi} \nabla \eta (\tU)\right)\,d \xi,\\
		&\mathcal{I}_5:=\int_{\mathbb{R}_+} (\nabla \eta )(U|\tU)\rd_{\xi}\left( M(\tU)\rd_{\xi} \nabla \eta(\tU)\right)\,d\xi.				
		\end{align*}
		Each of term above can be written as
		\begin{align}
		\begin{aligned}\label{est_v_h1_inflow}
		\mathcal{Y}(U)&=\int_{\mathbb{R}_+}\widetilde{h}_{\xi}(h-\widetilde{h})\,d\xi  -\int_{\mathbb{R}_+}p'(\widetilde{v})\widetilde{v}_{\xi}(v-\widetilde{v})\,d\xi,\\
		\mathcal{I}_1&=-\int_{\mathbb{R}_+}\partial_{\xi}\left((p(v)-p(\widetilde{v}))(h-\widetilde{h})-\sigma_-\eta(U|\widetilde{U})\right)\,d\xi \\
		& =\left[(p(v)-p(\widetilde{v}))(h-\widetilde{h})-{\sigma_-\over 2}|h-\widetilde{h}|^2 -\sigma_-Q(v|\tv)\right]_{\xi=0}, \\
		\mathcal{I}_2&=-\int_{\mathbb{R}_+}\widetilde{h}_{\xi} p(v|\widetilde{v})\,d\xi, \\
		\mathcal{I}_3&=\int_{\mathbb{R}_+}\left(p(v)-p(\widetilde{v})\right)\partial_{\xi}\left(\frac{1}{\gamma p(v)}\partial_{\xi}\left(p(v)-p(\widetilde{v})\right)\right)\,d\xi\\
		&=\left[-\frac{1}{\gamma p(v)}(p(v)-p(\widetilde{v}))\partial_{\xi}(p(v)-p(\widetilde{v}))\right]_{\xi=0}-\int_{\mathbb{R}_+}\frac{1}{\gamma p(v)} |(p(v)-p(\widetilde{v}))_{\xi}|^2\,d\xi, \\
		\mathcal{I}_4&=\int_{\mathbb{R}_+}\left(p(v)-p(\widetilde{v})\right)\partial_{\xi}\left(\frac{p(\widetilde{v})-p(v)}{\gamma p(v)p(\widetilde{v})}\partial_{\xi}p(\widetilde{v})\right)\,d\xi \\
		& =\left[\left(p(v)-p(\widetilde{v})\right)^2\frac{\partial_{\xi}p(\widetilde{v})}{\gamma p(v)p(\widetilde{v})}\right]_{\xi=0}-\int_{\mathbb{R}_+}\partial_{\xi}\left(p(v)-p(\widetilde{v})\right) \frac{p(\widetilde{v})-p(v)}{\gamma p(v)p(\widetilde{v})}\partial_{\xi}p(\widetilde{v})\,d\xi, \\
		\mathcal{I}_5&=-\int_{\mathbb{R}_+}p(v|\widetilde{v})(\ln \widetilde{v})_{\xi\xi}\,d\xi.
		\end{aligned}
		\end{align}
		
		Compared to the proof of Lemma \ref{lem:v-H1}, the only additional terms are included in $\mathcal{I}_{1}$, namely,
		\begin{equation}\label{new_terms}
		\left[-\frac{\sigma_-}{2}|h-\widetilde{h}|^2-\sigma_-Q(v|\tv)\right]_{\xi=0}.
		\end{equation}
		However, we use the definition of $h$ and $\widetilde{h}$ to observe that
		
		\begin{equation}\label{H3}
		\begin{aligned}
		\int_{0}^{t}\left[|h-\widetilde{h}|^2\right]_{\xi=0}\,ds & \leq C\int_{0}^{t} \left[|u-\widetilde{u}|^{2}  +|\widetilde{v}_{\xi}(v-\widetilde{v})|^{2}\right]_{\xi=0}\,ds +C\int_{0}^{t} [|(v-\widetilde{v})_{\xi}|^{2}]_{\xi=0}\,ds \\
		& \leq C\varepsilon e^{-C\delta\beta} + C\varepsilon\delta^{3}e^{-C\delta\beta} + C\int_{0}^{t} \left[|(v-\widetilde{v})_{\xi}|^{2}\right]_{\xi=0}\,ds.  
		\end{aligned}
		\end{equation}
		From the equations, we obtain the equation for the perturbation $v-\tv$ as
		\begin{align*}
		(v-\widetilde{v})_t -\sigma_-(v-\widetilde{v})_{\xi}   -(u-\widetilde{u})_{\xi} = \sX\widetilde{v}_{\xi},
		\end{align*}
		from which we derive
		\begin{equation}\label{1stbvv}
		\begin{aligned}
		&C\int_{0}^{t} \left[|(v-\widetilde{v})_{\xi}|^{2}\right]_{\xi=0}\,ds   \\
		&\leq {C}\int_{0}^{t} \left(\left[(v-\widetilde{v})^{2}_{t}\right]_{\xi=0} +  |\sX(s)|^{2}\left[|\widetilde{v}_{\xi}|^{2}\right]_{\xi=0}+ \left[(u-\widetilde{u})^{2}_{\xi}\right]_{\xi=0}\right)\,ds \\
		&\leq Ce^{-C\delta\beta} + C\varepsilon\delta^{3}e^{-C\delta\beta} +C\int_{0}^{t}\|(u-\widetilde{u})_{\xi}\|^{2}_{L^{\infty}(\mathbb{R}_+)}ds  \\
		& \leq Ce^{-C\delta\beta}  + C\int_{0}^{t} \|(u-\widetilde{u})_{\xi}\|_{L^2(\mathbb{R}_+)}\|(u-\widetilde{u})_{\xi\xi}\|_{L^{2}(\mathbb{R}_+)}\,ds \\
        & \leq Ce^{-C\delta\beta}  + C\int_{0}^{t} \|(u-\widetilde{u})_{\xi}\|^{2}_{L^{2}(\mathbb{R}_+)}\,ds + \frac{1}{100}\int_{0}^{t} \|(u-\widetilde{u})_{\xi\xi}\|^{2}_{L^{2}(\mathbb{R}_+)}\,ds \\
		&	= Ce^{-C\delta\beta} +  C\int_{0}^{t} D_{u_1}\,ds+ \frac{1}{100}\int_{0}^{t}D_{u_2}\,ds.
		\end{aligned}
		\end{equation}
		 We now substitute \eqref{1stbvv} into \eqref{H3} to derive 
		\begin{align*}
		\int_{0}^{t} \left[|h-\widetilde{h}|^2\right]_{\xi=0}\,ds  \leq Ce^{-C\delta\beta} +  C\int_{0}^{t} D_{u_1}\,ds + \frac{1}{100}\int_{0}^{t}D_{u_2}\,ds.
		\end{align*}
		On the other hand, it is easy to observe that
		\begin{align*}
		\left|\int_{0}^{t} [Q(v|\tv)]_{\xi=0} \,ds\right|\le C\lVert v-\widetilde{v} \rVert_{L^\infty(0, T; L^2(\mathbb{R}_+))} \int_{0}^{t} |\widetilde{v}(s,0)-v_-|\,ds \le C\e e^{-C\delta\beta}.
		\end{align*}
		Therefore, the above computation shows that the time integration of the new terms in \eqref{new_terms} can be bounded by
		\[Ce^{-C\delta\beta}+C\int_0^tD_{u_1}\,ds+\frac{1}{100}\int_0^tD_{u_2}\,ds.\]
		Since the other terms in \eqref{est_v_h1_inflow} are the same as \eqref{est_v_h1}, we deduce that the desired estimate can be obtained by using exactly the same estimate as in Lemma \ref{lem:v-H1}.
	\end{proof}
	Finally, we close the estimate on the $H^1$-perturbation in the following lemma.
	\begin{lemma}\label{estimateux}
		Under the hypothesis of Proposition \ref{prop:main}, there exists a positive constant $C$ such that for all $t\in[0,T]$,
		\begin{align*}
		&\lVert (v-\widetilde{v})(t, \cdot) \rVert^2_{H^1(\mathbb{R}_+)}+\lVert (u-\widetilde{u})(t, \cdot) \rVert^2_{H^1(\mathbb{R}_+)} +\delta \int_{0}^{t}|\sX|^2\,ds \\
		&+\int_{0}^{t}\left(G_1+G^S+D_{v_1}+D_{u_1}+D_{u_2}\right)\,ds\\   
		& \le  C\left(\lVert (v-\widetilde{v})(0, \cdot) \rVert^2_{H^1(\mathbb{R}_+)}+\lVert (u-\widetilde{u})(0, \cdot) \rVert^2_{H^1(\mathbb{R}_+)}\right)
		+ Ce^{-C\delta\beta}.
		\end{align*}
	\end{lemma}

	\begin{proof}
		We first observe that the perturbation $u-\tu$ satisfies
		\[(u-\widetilde{u})_t=\sigma_-(u-\widetilde{u})_{\xi}-(p(v)-p(\widetilde{v}))_{\xi}+\left(\frac{u_{\xi}}{v}-\frac{\widetilde{u}_{\xi}}{\widetilde{v}}\right)_{\xi}+\sX(t)\widetilde{u}_{\xi}.\]
		We multiply the above equation by $-(u-\widetilde{u})_{\xi\xi}$ and integrate it over $\left[0, t\right]\times \mathbb{R}_+$ to get
		\begin{align}
		\begin{aligned} \label{est_u_h1}
		-\int_{0}^{t}&\int_{\R_+}(u-\widetilde{u})_t(u-\widetilde{u})_{\xi\xi}\,d\xi\,ds\\
		&= \int_{0}^{t}\int_{\R_+} (p(v)-p(\widetilde{v}))_{\xi}(u-\widetilde{u})_{\xi\xi}\,d\xi\,ds -\int_{0}^{t}\int_{\R_+}\left(\frac{u_{\xi}}{v}-\frac{\widetilde{u}_{\xi}}{\widetilde{v}}\right)_{\xi}(u-\widetilde{u})_{\xi\xi}\,d\xi\,ds  \\
		&\quad -\int_{0}^{t}\int_{\R_+}\sX(s)\widetilde{u}_{\xi}(u-\widetilde{u})_{\xi\xi}\,d\xi \,ds -\sigma_{-}\int_{0}^{t}\int_{\R_+} (u-\widetilde{u})_{\xi}(u-\widetilde{u})_{\xi\xi}\,d\xi \,ds \\
		& =:J_1+J_2+J_3+J_4.
		\end{aligned}
		\end{align}
		
		Applying integration-by-parts for the left-hand side of \eqref{est_u_h1}, we have
		\begin{align*}
		&-\int_{0}^{t}\int_{0}^{\infty}(u-\widetilde{u})_{t}(u-\widetilde{u})_{\xi\xi}\,d\xi \,ds  \\
		&=\int_{0}^{t}\left[(u-\widetilde{u})_t(u-\widetilde{u})_{\xi} \right]_{\xi=0}\,ds+\frac{1}{2}\lVert (u-\widetilde{u})_{\xi}(t,\cdot) \rVert^2_{L^2(\mathbb{R}_+)}-\frac{1}{2}\lVert (u-\widetilde{u})_{\xi}(0,\cdot) \rVert^2_{L^2(\mathbb{R}_+)}.
		\end{align*}
		Thus, we rewrite \eqref{est_u_h1} as 
		\begin{equation}\label{est_u_h1_2}
		\frac{1}{2}\lVert (u-\widetilde{u})_{\xi}(t,\cdot) \rVert^2_{L^2(\mathbb{R}_+)}=\frac{1}{2}\lVert (u-\widetilde{u})_{\xi}(0,\cdot) \rVert^2_{L^2(\mathbb{R}_+)}+J_1+J_2+J_3 + J_4+J_5,
		\end{equation}
		where $J_5:=-\int_{0}^{t}[(u-\widetilde{u})_{t}(u-\widetilde{u})_{\xi}]_{\xi=0}\,ds$. Observe that the terms $J_1$, $J_2$, and $J_3$ are exactly the same terms as in \eqref{F1F2F3}, except that the variable is changed to $\xi$. Therefore, these terms can be treated exactly in the same manner as before. Thus, we focus on estimating $J_4$ and $J_5$. \\
		
		\noindent $\bullet$ (Estimate of $J_4$): By virtue of Cauchy-Schwarz inequality and Young's inequality, we have 
		\begin{align*}
		|J_4| &\leq C\int_{0}^{t}\|(u-\widetilde{u})_{\xi}\|_{L^{2}(\mathbb{R}_+)}\|(u-\widetilde{u})_{\xi\xi}\|_{L^{2}(\mathbb{R}_+)}\,ds \\
		& \leq   C\int_{0}^{t} D_{u_1}\,ds + \frac{1}{100}\int_{0}^{t} D_{u_2}\,ds.
		\end{align*}
		
		\noindent $\bullet$ (Estimate of $J_5$): Using the equation
		\[(v-\tv)_t-\sigma_-(v-\tv)_\xi - (u-\tu)_\xi = \sX\tv_\xi,\]
		we split $J_5$ as
		\begin{align*}
		J_5&=-\int_{0}^{t}\left[(u-\widetilde{u})_{t}(v-\widetilde{v})_{t}\right]_{\xi=0}\,ds +\sigma_-\int_{0}^{t}\left[(u-\widetilde{u})_{t}(v-\widetilde{v})_{\xi}\right]_{\xi=0}\,ds+\int_{0}^{t}\sX(s)\left[(u-\widetilde{u})_t\widetilde{v}_{\xi}\right]_{\xi=0}\,ds\\
		&=: J_{51} + J_{52} + J_{53}.
		\end{align*}
		For $J_{51}$, using integration-by-parts in time variable, we have
		\begin{align*}
		|J_{51}|&\le\left|((u-\widetilde{u})_{t}(v-\widetilde{v}))(t,0)-((u-\widetilde{u})_{t}(v-\widetilde{v}))(0,0)\right|+ \left|\int_{0}^{t}\left[(u-\widetilde{u})_{tt}(v-\widetilde{v})\right]_{\xi=0}\,ds\right| \\
		&\le C\|v-\tv\|_{L^\infty(0, T ;L^2(\mathbb{R}_+))}\left(|\tu_\xi(t,0)|+|\tu_\xi(0,0)|+\int_0^t|\tu_{\xi\xi}(s,0)|\,ds\right)\\
		&\leq C\varepsilon \delta^{2} e^{-C\delta\beta}.
		\end{align*}
		By Young's inequality and the same method as done in \eqref{1stbvv}, we estimate $F_{52}$ as 
		\begin{align*}
		|J_{52}|&\leq C\int_{0}^{t}\left[(u-\widetilde{u})^{2}_{t}\right]_{\xi=0}\,ds  + C\int_{0}^{t}\left[(v-\widetilde{v})^{2}\right]_{\xi=0}\,ds \\
		& \leq Ce^{-C\delta\beta} + C\varepsilon\delta^{3}e^{-C\delta\beta} + C\int_{0}^{t} D_{u_1}\,ds + \frac{1}{100}\int_{0}^{t} D_{u_2}\,ds.
		\end{align*}
		Similarly, we can obtain $|J_{53}|\leq Ce^{-C\delta\beta}$ by using the fact that $|\sX| \leq C\varepsilon$. Substituting the above estimates into \eqref{est_u_h1_2}, there exists a positive constant $C_4$ such that 
		\begin{align*}
		&\lVert (u-\widetilde{u})_{\xi}(t,\cdot) \rVert^2_{L^2(\mathbb{R}_+)}+\int_{0}^{t}D_{u_2}\,ds \\
		&\le C\lVert (u-\widetilde{u})_{\xi}(0,\cdot) \rVert^2_{L^2(\mathbb{R}_+)}+Ce^{-C\delta\beta}+C_4\int_{0}^{t} (\delta|\sX|^2+G_1+G^S+D_{u_1}+D_{v_1}) \,ds.
		\end{align*}
		Multiplying the above inequality by the constant $\frac{1}{2\max\{1,C_4\}}$, and then adding the result to \eqref{firstorderv}, we have 
		\begin{align*}
		\begin{aligned} 
		&\lVert (v-\widetilde{v})(t, \cdot) \rVert^2_{H^1(\mathbb{R}_+)}+\lVert (u-\widetilde{u})(t, \cdot) \rVert^2_{H^1(\mathbb{R}_+)} +\delta \int_{0}^{t}|\sX|^2 \,ds  \\
		&\quad +\int_{0}^{t}\left(G_1+G^S+D_{v_1}+D_{u_1}+D_{u_2}\right)\,ds\\   
		& \le  C\left(\lVert (v-\widetilde{v})(0, \cdot) \rVert^2_{H^1(\mathbb{R}_+)}+\lVert (u-\widetilde{u})(0, \cdot) \rVert^2_{H^1(\mathbb{R}_+)}\right)
		+ Ce^{-C\delta\beta}.
		\end{aligned}
		\end{align*}
		This implies the desired result in Lemma \ref{estimateux}, and completes the proof.
	\end{proof}

\begin{appendix}
\setcounter{equation}{0}
\section{Global existence}
	Based on the a priori estimates, we present the global existence of \eqref{eq:NS} and \eqref{eq:NS_inflow_2}. Since the proofs are almost identical, we focus on the impermeable wall problem \eqref{eq:NS}. We choose smooth monotone functions $\underline{v}$ and $\underline{u}$ defined on $\bbr_+$ such that
	\begin{align*}
	&\|\underline{v}-v_+\|_{L^2(\beta,\infty)}+\|\underline{v}-v_-\|_{L^2(0,\beta)}\\
	&\quad+\|\underline{u}-u_+\|_{L^2(\beta,\infty)}+\|u-u_-\|_{L^2(0,\beta)}+\|\underline{v}_x\|_{L^2(\R_+)}+\|\underline{u}_x\|_{L^2(\R_+)}<\underline{C}\delta,
	\end{align*}
	for some constant $\underline{C}$. Then, noticing $(\tv^{\bX,\beta},\tu^{\bX,\beta})(0,x)=(\tv,\tu)(x-\beta)$ and using \eqref{shock_prop}, we have
	\begin{align*}
		&\|\underline{v}(\cdot)-\tv^{\bX,\beta}(0,\cdot)\|_{H^1(\R_+)}+\|\underline{u}(\cdot)-\tu^{\bX,\beta}(0,\cdot)\|_{H^1(\R_+)}\\
		&\le\|\underline{v}-v_+\|_{L^2(\beta,\infty)}+\|\underline{v}-v_-\|_{L^2(0,\beta)}+\|\underline{u}-u_+\|_{L^2(\beta,\infty)}+\|\underline{u}-u_-\|_{L^2(0,\beta)}\\
		&\quad +\|\tv^{\bX,\beta}(0,\cdot)-v_+\|_{L^2(\beta,\infty)}+\|\tv^{\bX,\beta}(0,\cdot)-v_-\|_{L^2(0,\beta)}\\
		&\quad +\|\tu^{\bX,\beta}(0,\cdot)-u_+\|_{L^2(\beta,\infty)}+\|\tu^{\bX,\beta}(0,\cdot)-u_-\|_{L^2(0,\beta)}\\
		&\quad +\|\underline{v}_x\|_{L^2(\R_+)}+\|\pa_x\tv^{\bX,\beta}(0,\cdot)\|_{L^2(\R_+)} +\|\underline{u}_x\|_{L^2(\R_+)}+\|\pa_x\tu^{\bX,\beta}(0,\cdot)\|_{L^2(\R_+)}\\
		&\le C_1\sqrt{\delta},
	\end{align*}
	for another positive constant $C_1$. We now define a positive constant $\e_0$ and $\e_*$ as
	\[\e_0=\e_*-\underline{C}\delta,\quad \e_*:={\e\over 2(C_0 +1)}-C_1\sqrt{\delta}-e^{-C\delta\beta} ,\]
	where the constants $\e$ and $C_0$ are the same constants in \eqref{apriori_small} and \eqref{apriori_impermeable}. Here, $\e_0$ and $\e_*$ can be chosen as positive constant and $\e_*<\frac{\e}{2}$, thanks to the smallness of $\delta$ and $e^{-\beta}$. Now, consider any initial data $(v_0,u_0)$ that satisfies \eqref{initial_condi} in Theorem \ref{thm:impermeable}, that is,
	\begin{align*}
	&\|v_0-v_+\|_{L^2(\beta,\infty)}+\|v_0-v_-\|_{L^2(0,\beta)}\\
	&\quad +\|u_0-u_+\|_{L^2(\beta,\infty)}+\|u_0-u_-\|_{L^2(0,\beta)}+\|v_{0x}\|_{L^2(\R_+)}+\|u_{0x}\|_{L^2(\R_+)}<\e_0.
	\end{align*}
	Then, this implies
	\begin{align*}
		&\|v_0-\underline{v}\|_{H^1(\R_+)}+\|u_0-\underline{u}\|_{H^1(\R_+)}\\
		&\le \|v_0-v_+\|_{L^2(\beta,\infty)}+\|u_0-u_+\|_{L^2(\beta,\infty)}+\|\underline{v}-v_+\|_{L^2(\beta,\infty)}+\|\underline{u}-u_+\|_{L^2(\beta,\infty)}\\
		&\quad + \|v_0-v_-\|_{L^2(0,\beta)}+\|u_0-u_-\|_{L^2(0,\beta)}+\|\underline{v}-v_-\|_{L^2(0,\infty)}+\|\underline{u}-u_-\|_{L^2(0,\beta)}\\
		&\quad +\|v_{0x}\|_{L^2(\R_+)}+\|u_{0x}\|_{L^2(\R_+)}+\|\underline{v}_x\|_{L^2(\R_+)}+\|\underline{u}_x\|_{L^2(\R_+)}\\
		&\le \e_0+\underline{C}\delta = \e_*.
	\end{align*}
	In particular, Sobolev inequality implies $\|v_0-\underline{v}\|_{L^\infty(\R_+)}<C\e_*$ and therefore, for small enough $\e_*$, 
	\[\frac{v_-}{2}<v_0(x)<2v_+,\quad x\in\R_+.\]
	Therefore, the local existence theorem implies that there exists $T_0>0$ such that the system \eqref{eq:NS} admits a unique solution $(v,u)$ on $[0,T_0]$ satisfying
	\begin{equation}\label{vvunder}
		\|v-\underline{v}\|_{L^\infty(0,T_0;H^1(\R_+))}+\|u-\underline{u}\|_{L^\infty(0,T_0;H^1(\R_+))}\le \frac{\e}{2}
	\end{equation}
	and
	\[\frac{v_-}{3}<v(t,x)<3v_+,\quad t>0,\quad x\in\R_+.\]
	Then,
	\begin{align*}
		&\|\underline{v}-\tv^{\bX,\beta}(t,\cdot)\|_{H^1(\R_+)}+\|\underline{u}-\tu^{\bX,\beta}(t,\cdot)\|_{H^1(\R_+)}\\
		&\le \|\underline{v}-v_+\|_{L^2(\beta,\infty)}+\|\underline{u}-u_+\|_{L^2(\beta,\infty)}+\|\tv^{\bX,\beta}(t,\cdot)-v_+\|_{L^2(\beta,\infty)}+\|\tu^{\bX,\beta}(t,\cdot)-u_+\|_{L^2(\beta,\infty)}\\
		&\quad + \|\underline{v}-v_-\|_{L^2(0,\beta)}+\|\underline{u}-u_-\|_{L^2(0,\beta)}+\|\tv^{\bX,\beta}(t,\cdot)-v_-\|_{L^2(0,\beta)}+\|\tu^{\bX,\beta}(t,\cdot)-u_-\|_{L^2(0,\beta)}\\
		&\quad +\|\pa_x\tv^{\bX,\beta}(t,\cdot)\|_{L^2(\R_+)}+\|\pa_x\tu^{\bX,\beta}(t,\cdot)\|_{L^2(\R_+)}+\|\underline{v}_x\|_{L^2(\R_+)}+\|\underline{u}_x\|_{L^2(\R_+)}\\
		&\le C\sqrt{\delta}(1+\sqrt{|X(t)|})\le C\sqrt{\delta}(1+\sqrt{t}).
	\end{align*}
	Therefore, if we choose $\delta$ and $T_1\in(0,T_0)$ small enough so that $C\sqrt{\delta}(1+\sqrt{T_1})<\frac{\e}{2}$, we have
	\[\|\underline{v}-\tv^{\bX,\beta}(t,\cdot)\|_{L^\infty(0,T_1;H^1(\R_+))}+\|\underline{u}-\tu^{\bX,\beta}(t,\cdot)\|_{L^\infty(0,T_1;H^1(\R_+))}\le \frac{\e}{2}.\]
	Combining with the estimate \eqref{vvunder}, we obtain
	\[\|v-\tv^{\bX,\beta}(t,\cdot)\|_{L^\infty(0,T_1;H^1(\R_+))}+\|u-\tu^{\bX,\beta}(t,\cdot)\|_{L^\infty(0,T_1;H^1(\R_+))}\le \e.\]
	Moreover, since the shift $X$ is absolutely continuous, and
	\[v-\underline{v},u-\underline{u}\in C([0,T_1];H^1(\R_+)),\]
	we obtain $v-\tv,u-\tu\in C([0,T_1];H^1(\R_+))$. Therefore, we prove that \eqref{apriori_small} holds in the time interval $[0,T_1]$. We now show that the solution can be globally extended by using the standard continuation argument. To this end, suppose that the maximal existence time
	\[T_M:=\sup\left\{t>0~\Bigg|~\sup_{t\in[0,T]}\left(\|v-\tv^{\bX,\beta}\|_{H^1(\R_+)}+\|u-\tu^{\bX,\beta}\|_{H^1(\R_+)}\right)\le \e\right\}\]
	is finite. Then, we have
	\[\sup_{t\in[0,T_M]}\left(\|v-\tv^{\bX,\beta}\|_{H^1(\R_+)}+\|u-\tu^{\bX,\beta}\|_{H^1(\R_+)}\right)=\e.\]
	On the other hand, since
	\begin{align*}
	&\|v_0-\tv^{\bX,\beta}(0,\cdot)\|_{H^1(\R_+)}+\|u_0-\tu^{\bX,\beta}(0,\cdot)\|_{H^1(\R_+)}\\
	&=\|v_0-\underline{v}\|_{H^1(\R_+)}+\|u_0-\underline{u}\|_{H^1(\R_+)}+\|\underline{v}-\tv^{\bX,\beta}(0,\cdot)\|_{H^1(\R_+)}+\|\underline{u}-\tu^{\bX,\beta}(0,\cdot)\|_{H^1(\R_+)}\\
	&<\e_*+C_1\sqrt{\delta},
	\end{align*}
	the a priori estimate \eqref{apriori_impermeable} implies that
	\[\sup_{t\in[0,T_M]}\left(\|v-\tv^{\bX,\beta}\|_{H^1(\R_+)}+\|u-\tu^{\bX,\beta}\|_{H^1(\R_+)}\right)\le C_0(\e_*+C_1\sqrt{\delta})+C_0e^{-C\delta\beta}<\frac{\e}{2},\]
	which is a contradiction. Therefore, we obtain $T_M=+\infty$, and in particular, the a priori estimate holds for the whole time interval:
	\begin{align}
		\begin{aligned}\label{apriori_global}
			\sup_{t\ge0}&\left[\|v-\tv^{\bX,\beta}\|_{H^1(\R_+)}+\|u-\tu^{\bX,\beta}\|_{H^1(\R_+)}\right]\\
			&\quad +\sqrt{\int_0^\infty (\delta|\sX|^2+G_1+ G^S+D_{v_1}+D_{u_1}+D_{u_2})\,ds}\\
			&\le C_0\left(\|v_0-\tv^{\bX,\beta}(0,\cdot)\|_{H^1(\R_+)}+\|u_0-\tu^{\bX,\beta}(0,\cdot)\|_{H^1(\R_+)}\right)+C_0e^{-C\delta\beta},
		\end{aligned}
	\end{align}
	moreover,  for all $t>0$,
 \begin{equation} \label{bddx121}
|\dot{\bX}(t)|\le C_0 \left(\lVert (v-\widetilde{v}^{\bX,\beta})(t,\cdot) \rVert_{L^\infty(\mathbb{R}_+)}+\lVert (u-\widetilde{u}^{\bX,\beta})(t,\cdot) \rVert_{L^\infty(\mathbb{R}_+)}\right).
\end{equation}

	\section{Long-time behavior}
	We now show that the global-in-time estimate for the perturbation implies the desired long-time behaviors \eqref{long_time_impermeable} and \eqref{long_time_inflow} and thereby complete the proofs of Theorem \ref{thm:impermeable} and Theorem \ref{thm:inflow}. First, we define 
	\[g(t):=\|(v-\tv^{\bX,\beta})_x\|_{L^2(\R_+)}^2+\|(u-\tu^{\bX,\beta})_x\|_{L^2(\R_+)}^2,\]
	for the impermeable wall problem, and the same quantity is defined for inflow problem by substituting the variable $x$ by $\xi$. Once we show that $g\in W^{1,1}(\R_+)$, we have
	\[\lim_{t\to\infty}g(t)=0.\] 
	Then the interpolation inequality and the uniform bound on the $L^2$-perturbation in the a priori estimate imply that
	\begin{align*}
	\|v-\tv^{\bX,\beta}\|_{L^\infty(\R_+)}\le \|v-\tv^{\bX,\beta}\|_{L^2(\R_+)}^{1/2}\|(v-\tv^{\bX,\beta})_x\|_{L^2(\R_+)}^{1/2}\to 0,\quad\mbox{as}\quad t\to\infty,
	\end{align*}
	and the same estimate holds for $u-\tu^{\bX,\beta}$. Therefore, we only need to show that $g\in W^{1,1}(\R_+)$. Due to the different boundary conditions, we need to prove impermeable wall problem and inflow problem separately.
	
	\subsection{Impermeable wall problem}
	We first show that $g\in W^{1,1}(\R_+)$ for the impermeable wall condition problem.\\
	
	\noindent $\bullet$ ($g\in L^1(\R_+)$): We first observe that
	\begin{align*}
	(p(v)-p(\tv^{\bX,\beta}))_x &= p'(v)(v-\tv^{\bX,\beta})_x+\tv^{\bX,\beta}_x(p'(v)-p'(\tv^{\bX,\beta})),
	\end{align*}
	which yields
	\[|(v-\tv^{\bX,\beta})_x|\le C|(p(v)-p(\tv^{\bX,\beta}))_x|+C|\tv^{\bX,\beta}_x||v-\tv^{\bX,\beta}|.\]
	Therefore, using the relation $\sigma|\tv_\zeta|=|\tu_\zeta|$ and \eqref{apriori_global}, we obtain
	\[\int_0^\infty |g(t)|\,dt\le C\int_0^\infty \left(G^S+D_{v_1}+D_{u_1}\right)\,dt<+\infty.\]
	
	\noindent $\bullet$ ($g'\in L^1(\R_+)$): To show that $g'\in L^1(\R_+)$, we first note that the perturbations $v-\tv^{\bX,\beta}$ and $u-\tu^{\bX,\beta}$ satisfies
	
	\begin{align}
	\begin{aligned}\label{eq:perturb}
	&(v-\tv^{\bX,\beta})_t -(u-\tu^{\bX,\beta})_x=\sX\tv^{\bX,\beta}_x,\\
	&(u-\tu^{\bX,\beta})_t +(p(v)-p(\tv^{\bX,\beta}))_x=\left(\frac{u_x}{v}-\frac{\tu^{\bX,\beta}_x}{\tv^{\bX,\beta}}\right)_x+\sX\tu^{\bX,\beta}_x.
	\end{aligned}
	\end{align}
	Using the equations \eqref{eq:perturb} and integration-by-parts, we obtain
	
	\begin{align*}
	\int_0^\infty |g'(t)|\,dt&=\int_0^\infty 2\left|\int_{\R_+}(v-\tv^{\bX,\beta})_x(v-\tv^{\bX,\beta})_{tx}\,dx+\int_{\R_+}(u-\tu^{\bX,\beta})_x(u-\tu^{\bX,\beta})_{tx}\,dx\right|\,dt\\
	&\le 2\int_0^\infty \left|\int_{\R_+}(v-\tv^{\bX,\beta})_x\left[(u-\tu^{\bX,\beta})_x+\sX\tv^{\bX,\beta}_x\right]\,dx\right|\,dt\\
	&\quad +2\int_0^\infty \left|\int_{\R_+}(u-\tu^{\bX,\beta})_{xx}\left[-(p(v)-p(\tv^{\bX,\beta}))_x+\left(\frac{u_x}{v}-\frac{\tu^{\bX,\beta}_x}{\tv^{\bX,\beta}}\right)_x+\sX\tu^{\bX,\beta}_x\right]\,dx\right|\,dt\\
	&\quad +2\int_0^\infty \left|\left[(u-\tu^{\bX,\beta})_x(u-\tu^{\bX,\beta})_t\right]_{x=0}\right|\,dt\\
	&\le 2\int_0^\infty \left(\delta|\sX|^2+G^S+D_{v_1}+D_{u_1}+D_{u_2}\right)\,dt+2\int_0^\infty \int_{\R_+}\left|\left(\frac{u_x}{v}-\frac{\tu^{\bX,\beta}_x}{\tv^{\bX,\beta}}\right)_x\right|^2\,dxdt\\
	&\quad + 2\int_0^\infty \left|\left[(u-\tu^{\bX,\beta})_x(u-\tu^{\bX,\beta})_t\right]_{x=0}\right|\,dt.
	\end{align*}
	However, the second term on the right-hand side can be estimated by using the same argument as in \cite{KVW23}, which yields
	\[\int_0^\infty \int_{\R_+}\left|\left(\frac{u_x}{v}-\frac{\tu^{\bX,\beta}_x}{\tv^{\bX,\beta}}\right)_x\right|^2\,dxdt\le C\int_0^\infty \left(G^S+D_{v_1}+D_{u_1}+D_{u_2}\right)\,dt<+\infty.\]
	On the other hand, using $|\dot{\bX}(t)|\le C$ by \eqref{apriori_global} and \eqref{bddx12}, and the interpolation inequality, the last term is estimated as
	\begin{align*}
	\int_0^\infty& \left|\left[(u-\tu^{\bX,\beta})_x(u-\tu^{\bX,\beta})_t\right]_{x=0}\right|\,dt\\
	&\le \int_0^\infty |(u-\tu^{\bX,\beta})_x(t,0)\tu^{\bX,\beta}_x(t,0)||\sigma+\sX(t)|\,dt\le C\int_0^\infty  \|(u-\tu^{\bX,\beta})_x\|_{L^\infty(\R_+)}|\tu^{\bX,\beta}_x(t,0)|\,dt\\
	&\le C\int_0^\infty \|(u-\tu^{\bX,\beta})_x\|_{L^2(\R_+)}^{1/2}\|(u-\tu^{\bX,\beta})_{xx}\|_{L^2(\R_+)}^{1/2}|\tu^{\bX,\beta}_x(t,0)|\,dt\\
	&\le C\int_0^\infty\left(\|(u-\tu^{\bX,\beta})_x\|_{L^2(\R_+)}^2 +\|(u-\tu^{\bX,\beta})_{xx}\|_{L^2(\R_+)}^2+|\tu^{\bX,\beta}_x(t,0)|^2\right)\,dt\\
	&\le C\int_0^\infty (D_{u_1}+D_{u_2}+\delta^4e^{-C\delta t})\,dt<+\infty.
	\end{align*}
	To sum up, we obtain
	\[\int_0^\infty |g'(t)|\,dt<+\infty,\]
	which completes the proof of $g\in W^{1,1}(\R_+)$.
	
	\subsection{Inflow problem}
	Since the proof of $g\in L^1(\R_+)$ is the same as in the impermeable wall problem, we focus on showing $g'\in L^1(\R_+)$. For the inflow problem, the perturbation satisfies
	\begin{align*}
	&(v-\tv^{\bX,\beta})_t-\sigma_-(v-\tv^{\bX,\beta})_\xi-(u-\tu^{\bX,\beta})_\xi=\sX\tv^{\bX,\beta}_\xi,\\
	&(u-\tu^{\bX,\beta})_t-\sigma_-(u-\tu^{\bX,\beta})_\xi+(p(v)-p(\tv^{\bX,\beta}))_\xi=\left(\frac{u_\xi}{v}-\frac{\tu^{\bX,\beta}_\xi}{\tv^{\bX,\beta}}\right)_\xi+\sX\tu^{\bX,\beta}_\xi.
	\end{align*}
	Using these equations, we obtain
	\begin{align*}
		\int_0^\infty |g'(t)|\,dt&=\int_0^\infty 2\left|\int_{\R_+}(v-\tv^{\bX,\beta})_{\xi}(v-\tv^{\bX,\beta})_{t\xi}\,d\xi+\int_{\R_+}(u-\tu^{\bX,\beta})_{\xi}(u-\tu^{\bX,\beta})_{t\xi}\,d\xi\right|\,dt\\
		&\le 2\int_0^\infty \left|\int_{\R_+}(v-\tv^{\bX,\beta})_{\xi}\left[ \sigma_-(v-\tv^{\bX,\beta})_{\xi\xi} +(u-\tu^{\bX,\beta})_{\xi\xi}+\sX\tv^{\bX,\beta}_{\xi\xi}\right]\,d\xi\right|\,dt\\
		&\quad +2\int_0^\infty \left|\int_{\R_+}(u-\tu^{\bX,\beta})_{\xi\xi}\left[-(p(v)-p(\tv^{\bX,\beta}))_{\xi}+\left(\frac{u_{\xi}}{v}-\frac{\tu^{\bX,\beta}_{\xi}}{\tv^{\bX,\beta}}\right)_{\xi}+\sX\tu^{\bX,\beta}_{\xi}\right]\,d\xi\right|\,dt\\
		&\quad +2\int_0^\infty \left|\left[(u-\tu^{\bX,\beta})_{\xi}(u-\tu^{\bX,\beta})_t\right]_{\xi=0}\right|\,dt + |\sigma_-|\int_{0}^{\infty} \left|\left[(v-\tv^{\bX,\beta})^2_{\xi}\right]_{\xi=0}\right|\,dt  \\
		&\le 2\int_0^\infty \left(\delta|\dot{\bold{X}}|^2+\mathcal{G}^{S}+\bold{D}_{v_1}+\bold{D}_{u_1}+\bold{D}_{u_2}\right)\,dt+2\int_0^\infty \int_{\R_+}\left|\left(\frac{u_{\xi}}{v}-\frac{\tu^{\bX,\beta}_{\xi}}{\tv^{\bX,\beta}}\right)_{\xi}\right|^2\,d\xi dt\\
		&\quad + 2\int_0^\infty \left|\left[(u-\tu^{\bX,\beta})_{\xi}(u-\tu^{\bX,\beta})_t\right]_{\xi=0}\right|\,dt + |\sigma_-|\int_{0}^{\infty} \left|\left[(v-\tv^{\bX,\beta})^2_{\xi}\right]_{\xi=0}\right|\,dt.
	\end{align*}
	Notice that, the only difference compared to the impermeable wall problem is the last term:
	\[\int_{0}^{\infty} \left|\left[(v-\tv^{\bX,\beta})^2_{\xi}\right]_{\xi=0}\right|\,dt.\]
	However, using the equation for $v-\tv^{\bX,\beta}$, the boundary condition $v(t,0)=v_-$, and \eqref{shock_prop}, we have
	\begin{align*}
		&\int_{0}^{\infty} \left|\left[(v-\tv^{\bX,\beta})^2_{\xi}\right]_{\xi=0}\right|\,dt \\
		&\leq {C}\int_{0}^{\infty} \left(\left[(v-\tv^{\bX,\beta})^{2}_{t}\right]_{\xi=0} +  |\sX(s)|^{2}\left[\left|\tv^{\bX,\beta}_{\xi}\right|^{2}\right]_{\xi=0}+ \left[(u-\tu^{\bX,\beta})^{2}_{\xi}\right]_{\xi=0}\right)\,dt\\
		&\leq Ce^{-C\delta\beta} + C\delta^{3}e^{-C\delta\beta} +C\int_{0}^{\infty}\|(u-\tu^{\bX,\beta})_{\xi}\|^{2}_{L^{\infty}(\mathbb{R}+)}\,dt  \\
		& \leq Ce^{-C\delta\beta}  + C\int_{0}^{\infty} \|(u-\tu^{\bX,\beta})_{\xi}\|_{L^2(\mathbb{R}+)}\|(u-\tu^{\bX,\beta})_{\xi\xi}\|_{L^{2}(\mathbb{R}+)}\,dt \\
        &\le Ce^{-C\delta\beta} +  C\int_{0}^{\infty} D_{u_1}\,dt+ C\int_{0}^{\infty}D_{u_2}\,dt<+\infty.
	\end{align*}
	Therefore, we show that $g'\in L^1(\R_+)$ and this completes the proof of $g\in W^{1,1}(\R_+)$.\\
	
	Hence, it holds from \eqref{bddx121} that
	 \[
|\dot{\bX}(t)|\le C_0 \left(\lVert (v-\widetilde{v}^{\bX,\beta})(t,\cdot) \rVert_{L^\infty(\mathbb{R}_+)}+\lVert (u-\widetilde{u}^{\bX,\beta})(t,\cdot) \rVert_{L^\infty(\mathbb{R}_+)}\right) \to 0 \quad \mbox{as } t\to \infty.
\]
\end{appendix}

	\bibliographystyle{amsplain}
	\bibliography{reference} 

\providecommand{\bysame}{\leavevmode\hbox to3em{\hrulefill}\thinspace}
\providecommand{\MR}{\relax\ifhmode\unskip\space\fi MR }
\providecommand{\MRhref}[2]{%
  \href{http://www.ams.org/mathscinet-getitem?mr=#1}{#2}
}
\providecommand{\href}[2]{#2}
\begin{thebibliography}{10}

\bibitem{D96}
C.~M. Dafermos, \emph{Entropy and the stability of classical solutions of
  hyperbolic systems of conservation laws}, Recent mathematical methods in
  nonlinear wave propagation ({M}ontecatini {T}erme, 1994), Lecture Notes in
  Math., vol. 1640, Springer, Berlin, 1996, pp.~48--69.

\bibitem{D79}
R.~J. DiPerna, \emph{Uniqueness of solutions to hyperbolic conservation laws},
  Indiana Univ. Math. J. \textbf{28} (1979), no.~1, 137--188.

\bibitem{G86}
J.~Goodman, \emph{Nonlinear asymptotic stability of viscous shock profiles for
  conservation laws}, Arch. Rational Mech. Anal. \textbf{95} (1986), no.~4,
  325--344.

\bibitem{HKK23}
S.~Han, M.-J. Kang, and J.~Kim, \emph{Large-time behavior of composite waves of
  viscous shocks for the barotropic {N}avier-{S}tokes equations}, SIAM J. Math.
  Anal. \textbf{55} (2023), no.~5, 5526--5574.

\bibitem{HKKL_pre}
S.~Han, M.-J. Kang, J.~Kim, and H.~Lee, \emph{Long-time behavior towards
  viscous-dispersive shock for {N}avier-{S}tokes equations of {K}orteweg type},
  arXiv preprint.

\bibitem{HLS10}
F.~Huang, J.~Li, and X.~Shi, \emph{Asymptotic behavior of solutions to the full
  compressible {N}avier-{S}tokes equations in the half space}, Commun. Math.
  Sci. \textbf{8} (2010), 639--654.

\bibitem{HMS03}
F.~Huang, A.~Matsumura, and X.~Shi, \emph{Viscous shock wave and boundary layer
  solution to an inflow problem for compressible viscous gas}, Commun. Math.
  Phys. \textbf{239} (2003), 261--285.

\bibitem{HQ09}
F.~Huang and X.~Qin, \emph{Stability of boundary layer and rarefaction wave to
  an outflow problem for compressible {N}avier--{S}tokes equations under large
  perturbation}, J. Differential Equations \textbf{246} (2009), no.~10,
  4077--4096.

\bibitem{KV16}
M.-J. Kang and A.~F. Vasseur, \emph{Criteria on contractions for entropic
  discontinuities of systems of conservation laws}, Arch. Rational Mech. Anal.
  \textbf{222} (2016), no.~1, 343--391.

\bibitem{KV21}
\bysame, \emph{Contraction property for large perturbations of shocks of the
  barotropic {N}avier-{S}tokes system}, J. Eur. Math. Soc. \textbf{23} (2021),
  no.~2, 585--638.

\bibitem{KV-Inven}
\bysame, \emph{Uniqueness and stability of entropy shocks to the isentropic
  {E}uler system in a class of inviscid limits from a large family of
  {N}avier-{S}tokes systems}, Invent. Math. \textbf{224} (2021), no.~1,
  55--146.

\bibitem{KV-2shock}
\bysame, \emph{Well-posedness of the {R}iemann problem with two shocks for the
  isentropic {E}uler system in a class of vanishing physical viscosity limits},
  J. Differential Equations \textbf{338} (2022), 128--226.

\bibitem{KVW23}
M.-J. Kang, A.~F. Vasseur, and Y.~Wang, \emph{Time-asymptotic stability of
  composite waves of viscous shock and rarefaction for barotropic
  {N}avier-{S}tokes equations}, Adv. Math. \textbf{419} (2023), Paper No.
  108963, 66.

\bibitem{KVW-NSF}
\bysame, \emph{Time-asymptotic stability of generic riemann solutions for
  compressible {N}avier-{S}tokes-{F}ourier equations}, arXiv preprint
  arXiv:2306.05604 (2023).

\bibitem{KNZ03}
S.~Kawashima, S.~Nishibata, and P.~Zhu, \emph{Asymptotic stability of the
  stationary solution to the compressible {N}avier--{S}tokes equations in the
  half space}, Commun. Math. Phys. \textbf{240} (2003), no.~3, 483--500.

\bibitem{KZ08}
S.~Kawashima and P.~Zhu, \emph{Asymptotic stability of nonlinear wave for the
  compressible {N}avier--{S}tokes equations in the half space}, J. Differential
  Equations \textbf{244} (2008), no.~12, 3151--3179.

\bibitem{L85}
T.-P. Liu, \emph{Nonlinear stability of shock waves for viscous conservation
  laws}, Bull. Amer. Math. Soc. \textbf{12} (1985), no.~2, 233--236.

\bibitem{M01}
A.~Matsumura, \emph{Inflow and outflow problems in the half space for a
  one-dimensional isentropic model system of compressible viscous gas}, Methods
  Appl. Anal. \textbf{8} (2001), no.~4, 645--666.

\bibitem{MM99}
A.~Matsumura and M.~Mei, \emph{Convergence to travelling fronts of solutions of
  the {$p$}-system with viscosity in the presence of a boundary}, Arch.
  Rational Mech. Anal. \textbf{146} (1999), no.~1, 1--22.

\bibitem{MN85}
A.~Matsumura and K.~Nishihara, \emph{On the stability of travelling wave
  solutions of a one-dimensional model system for compressible viscous gas},
  Japan J. Appl. Math. \textbf{2} (1985), no.~1, 17--25.

\bibitem{MN86}
\bysame, \emph{Asymptotics toward the rarefaction waves of the solutions of a
  one-dimensional model system for compressible viscous gas}, Japan J. Appl.
  Math. \textbf{3} (1986), no.~1, 1--13.

\bibitem{MN99}
\bysame, \emph{Global asymptotics toward the rarefaction wave for solutions of
  viscous $p$-system with boundary effect}, Quart. Appl. Math. \textbf{58}
  (2000), no.~1, 69--83.

\bibitem{MN01}
\bysame, \emph{Large-time behaviors of solutions to an inflow problem in the
  half space for a one-dimensional system of compressible viscous gas}, Commun.
  Math. Phys. \textbf{222} (2001), 449--474.

\bibitem{QW09}
X.~Qin and Y.~Wang, \emph{Stability of wave patterns to the inflow problem of
  full compressible {N}avier--{S}tokes equations}, SIAM J. Math. Anal.
  \textbf{41} (2009), no.~5, 2057--2087.

\bibitem{QW11}
\bysame, \emph{Large-time behavior of solutions to the inflow problem of full
  compressible {N}avier--{S}tokes equations}, SIAM J. Math. Anal. \textbf{43}
  (2011), no.~1, 341--366.

\bibitem{SZ93}
A.~Szepessy and Z.~Xin, \emph{Nonlinear stability of viscous shock waves},
  Arch. Rational Mech. Anal. \textbf{122} (1993), 53--103.

\end{thebibliography}
	
\end{document}